\documentclass[a4paper,twoside,11pt,openright]{article}

\usepackage{bea_article}

\begin{document}

\title{Conformal invariance of isoradial dimer models \\
\&\\ the case of triangular quadri-tilings}
\author{B\'eatrice de Tili\`ere
\thanks{Supported by Swiss National Fund grant 47102009}\vspace{0.4cm}\\
       {\small Institut f\"ur Mathematik,
       Universit\"at Z\"urich,
       Winterthurerstrasse 190,
       CH-8057 Z\"urich.}\\
       {\small\texttt{beatrice.detiliere@math.unizh.ch}}
}
\date{}
\maketitle

\vspace{-1cm}
\begin{abstract}
We consider dimer models on graphs which are bipartite, periodic and
satisfy a geometric condition called {\em isoradiality}, defined in
\cite{Kenyon3}. We show that the scaling limit of the height function 
of any such dimer model is $1/\sqrt{\pi}$ times a Gaussian
free field. Triangular quadri-tilings were introduced in \cite{Bea}; 
they are dimer models on a family of isoradial graphs
arising form rhombus tilings. By means of two height functions, they
can be interpreted as random interfaces in dimension $2+2$. We show
that the scaling limit of
each of the two height functions is $1/\sqrt{\pi}$ times a Gaussian free field,
and that the two Gaussian free fields are independent.
\end{abstract}

\section{Introduction}\label{sec0}

\subsection{Height fluctuations for isoradial dimer models}\label{subsec01}

\subsubsection{Dimer models}\label{subsubsec011}

The setting for this paper is the {\em dimer model}. It is a
statistical mechanics model representing diatomic molecules
adsorbed on the surface of a crystal. An interesting feature of the
dimer model is that it is one of the very few statistical mechanics
models where exact and explicit results can be obtained, see
\cite{Kenyon6,Kenyon5} for an overview. Another very interesting
aspect is the alleged conformal invariance of its scaling limit,
which is already proved in the domino and $60^{\circ}$-rhombus cases
\cite{Kenyon1,Kenyon2,Kenyon4}. Theorem \ref{thm1} of this paper shows
this property for a wide class of dimer models containing the
above two cases.\jump
In order to give some insight, let us precisely define the setting.
The {\bf dimer model} is in bijection with a
 mathematical model called the {\bf $2$-tiling model}
representing discrete random interfaces. The {\em system} considered
for a $2$-tiling model is a planar graph $G$. {\em Configurations} of the
system, or {\bf $2$-tilings}, are coverings of $G$ with polygons
consisting of pairs of edge-adjacent faces of $G$, also
called {\bf $2$-tiles}, which leave no hole and don't overlap.
 The {\em system} of the corresponding dimer model is the dual
graph $G^*$ of $G$. {\em Configurations} of the dimer model are {\bf perfect matchings} of $G^*$,
that is set of edges covering every vertex exactly once. Perfect
matchings of $G^*$ determine $2$-tilings of $G$ as explained by the
following correspondence. Denote by $f^*$ the dual vertex of a face
$f$ of $G$, and consider an edge $f^* g^*$ of $G^*$. We say that the
$2$-tile of $G$ made of the adjacent faces $f$ and $g$ is the $2$-tile
{\bf corresponding} to the edge $f^* g^*$. Then $2$-tiles
corresponding to edges of a dimer configuration form a $2$-tiling of $G$.
Let us denote by $\M(G^*)$ the set of dimer configurations of $G^*$.\jump
As for all statistical mechanics models, dimer configurations are chosen
with respect to the {\bf Boltzmann measure} defined as follows. Suppose that
the graph $G^*$ is finite, and that a positive weight function $\nu$
is assigned to edges of $G^*$, then each dimer configuration
$M$ has an {\em energy}, $\E(M)=-\sum_{e\in M}\log \nu(e)$.
The probability of occurrence of the dimer configuration $M$ chosen with
respect to the Boltzmann measure $\mu^1$ is:
\begin{equation*}
\mu^1(M)=\frac{e^{-\E(M)}}{Z(G^*,\nu)}=\frac{\prod_{e\in
M}\nu(e)}{Z(G^*,\nu)},
\end{equation*}
where $Z(G^*,\nu)$ is the normalizing constant called the {\bf partition function}.
Using the bijection between dimer configurations and $2$-tilings,
$\nu$ can be seen as weighting $2$-tiles, and $\mu^1$ as a
measure on $2$-tilings of $G$. When the graph $G^*$ is infinite, a
{\bf Gibbs measure} is defined to be a probability measure on $\M(G^*)$ with the
following property: if the matching in an annular region is fixed,
then matchings inside and outside of the annulus are independent,
moreover the probability of any interior matching $M$ is proportional to
$\prod_{e\in M}\nu(e)$. From now on, let us assume that the graph $G$
satisfies condition $(*)$
below:\jump
\begin{tabular}{p{0.6cm}|p{11.9cm}}
\vspace{0.22cm}$(*)$&The graph $G$ is infinite, planar,
and simple ($G$ has no loops and no multiple edges); its vertices
are of degree $\geq 3$. $G$ is {\em simply connected}, i.e. it is
the one-skeleton of a simply connected union of faces; and it is made 
of finitely many different faces, up to isometry.
\end{tabular}

\subsubsection{Isoradial dimer models}\label{subsubsec012}

This paper actually proves conformal invariance of the scaling limit
for a sub-family of all dimer models called {\em isoradial dimer
models}, introduced by Kenyon in \cite{Kenyon3}. Much attention has
lately been given to isoradial dimer models because of a
surprising feature: many statistical mechanics quantities can be
computed in terms of the {\em local} geometry of the graph. This fact
was conjectured in \cite{Kenyon3}, and proved in \cite{Bea1}. The
motivation for their study is further enhanced by the fact that the
yet classical domino and $60^{\circ}$-rhombus tiling models are examples of
isoradial dimer models. Last but not least their understanding
allows us to apprehend a random interface model in dimension $2+2$
called the {\em triangular quadri-tiling model} introduced in
\cite{Bea}, see Section \ref{subsubsec021}.\jump
Let us now define isoradial dimer models. Speaking in the terminology of $2$-tilings,
{\bf isoradial $2$-tiling models} are defined on graphs $G$ satisfying
a geometric condition called {\bf isoradiality}: all faces of an isoradial graph
are inscribable in a circle, and all circumcircles have the same
radius, moreover all circumcenters of the faces are contained in the
closure of the faces. The energy of configurations is determined
by a specific weight function called the {\em critical weight
function}, see Section \ref{subsec11} for definition.
Note that if $G$ is an isoradial graph, an isoradial embedding of the
dual graph $G^*$ is obtained by sending dual vertices to the center of
the corresponding faces. Hence, the corresponding dimer model is
called an {\bf isoradial dimer model}.

\subsubsection{Height functions}\label{subsubsec013}

Let $G$ be an isoradial graph, whose dual graph $G^*$ is bipartite.
Then, by means of the {\em height function}, $2$-tilings of $G$ can be
interpreted as random discrete $2$-dimensional surfaces in a $3$-dimensional
space that are projected orthogonally to the plane. In physics
terminology, one speaks of random interfaces in dimension $2+1$. The height
function, denoted by $h$, is an $\RR$-valued function on the vertices
 of every $2$-tiling of $G$, and is defined in Section \ref{sec2}.

\subsubsection{Gaussian free field in the plane}\label{subsubsec014}

The Gaussian free field in the plane is defined in Section \ref{sec3}. It is
a random distribution which assigns to functions
$\vphi_1,\ldots,\vphi_k\in\K$ (the set of compactly supported smooth
functions of $\RR^2$, which have mean $0$), a real Gaussian random vector
$(F\vphi_1,\ldots,F\vphi_k)$ whose covariance function is given by
\begin{equation*}
\EE(F\vphi_i
F\vphi_j)=\int_{\RR^2}\int_{\RR^2}g(x,y)\vphi_i(x)\vphi_j(y)dx\,dy,
\end{equation*}
where $g(x,y)=-\frac{1}{2\pi}\log|x-y|$ is the Green function of the
plane (defined up to an additive constant). The Gaussian
free field is conformally invariant \cite{Kenyon2}.

\subsubsection{Statement of result}\label{subsubsec015}

Let $G$ be an isoradial graph, whose dual graph $G^*$ is
bipartite. Suppose moreover that $G^*$ is doubly periodic, i.e. that
the graph $G^*$ and its vertex-coloring are periodic. Then by 
Sheffield's theorem \cite{Sheffield0}, there
exists a two-parameter family of translation invariant, ergodic Gibbs
measures; let us denote by $\mu$ the unique measure
which has minimal free energy per fundamental domain. From now on, we
assume that dimer configurations of $G^*$ are chosen with respect
to the measure $\mu$.\jump
Let us multiply the edge-lengths of the graph $G$ by $\eps>0$, this
yields a new graph $G^\eps$. Let $h^\eps$ be the unnormalized height
function on $2$-tilings of $G^\eps$. An important issue in the study
of the dimer model is the understanding of the fluctuations of $h^\eps$, as
the mesh $\eps$ tends to $0$. This question is answered by Theorem
\ref{thm1} below. Define:
\begin{equation*}
\begin{array}{lccl}
H^\eps : & \K & \rightarrow & \RR\\
 &\vphi& \longmapsto &\displaystyle H^\eps \vphi =
\eps^2 \sum_{v \in V(\Geps)}a(v^*)\vphi(v)h^\eps(v),
\end{array}
\end{equation*}
where $V(\Geps)$ denotes the set of vertices of the graph $\Geps$, and
$a(v^*)$ is the area in $G^*$ of the dual face $v^*$ of a vertex $v$.
\begin{thm}\label{thm1}
Consider a graph $G$ satisfying the above assumptions, then
$H^\eps$ converges weakly in distribution to
$\frac{1}{\sqrt{\pi}}$ times a Gaussian free field, that is for
every $\vphi_1,\ldots,\vphi_k \in \K$, $(H^\eps
\vphi_1,\ldots,H^\eps \vphi_k)$ converges in law (as $\eps
\rightarrow 0$) to $\frac{1}{\sqrt{\pi}}(F \vphi_1,\ldots,F
\vphi_k)$, where $F$ is a Gaussian free field.
\end{thm}
$\bullet$ As a direct consequence of Theorem \ref{thm1}, we obtain
convergence of the height function of domino and $60^{\circ}$-rhombus
tilings chosen with respect to the uniform measure to a Gaussian free field. Note that this
result is slightly different than those of
\cite{Kenyon1,Kenyon2,Kenyon4} since we work on the whole plane, and
not on simply connected regions.\jump
$\bullet$ The method for proving Theorem \ref{thm1} is
essentially that of \cite{Kenyon1}, except Lemma \ref{lem59} which is
new. Nevertheless, since we work with a general isoradial graph
(and not the square lattice), each step is adapted in a non
trivial way.

\subsection{The case of triangular quadri-tilings}\label{subsec02}

\subsubsection{Triangular quadri-tiling model}\label{subsubsec021}

\begin{figure}[ht]
\begin{center}
\includegraphics[height=10cm]{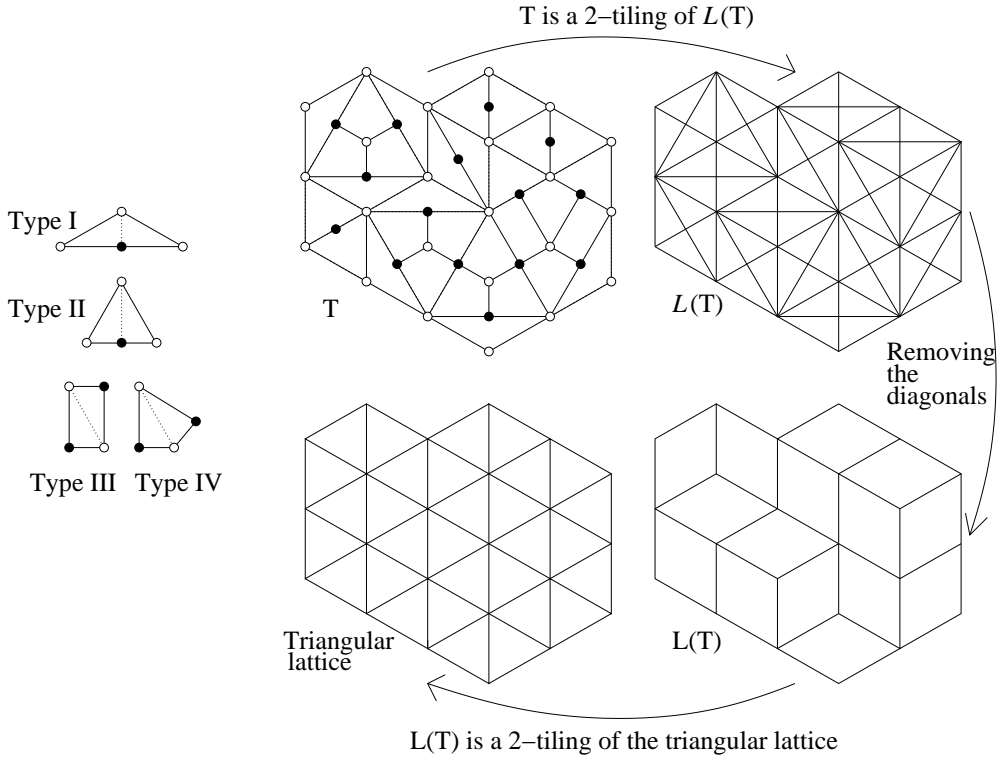}
\end{center}
\caption{Four type of quadri-tiles (left). Triangular quadri-tilings
correspond to two superposed dimer models (right).} \label{fig2}
\end{figure}

An exciting consequence of Theorem \ref{thm1} is that it allows us to
understand height fluctuations in the case of a random interface model in
dimension $2+2$, called the {\em triangular quadri-tiling model}. 
It is the first time this type of result can be obtained on such a
model.\jump
Let us start by defining triangular quadri-tilings. Consider the set of right
triangles whose hypotenuses have length $1$, and whose interior angle
is $\pi/3$. Color the vertex at the right angle black, and the other
two vertices white. A {\bf quadri-tile} is a quadrilateral obtained
from two such triangles in two different ways: either glue them
along a leg of the same length matching the black (white) vertex
to the black (white) one, or glue them along the hypotenuse. There
are four types of quadri-tiles classified as I, II, III, IV,
each of which has four vertices, see Figure \ref{fig2} (left). A {\bf
triangular quadri-tiling} of the
plane is an edge-to-edge tiling of the plane by quadri-tiles that
respects the coloring of the vertices, see $T$ of Figure \ref{fig2}
for an example. Let $\Q$ denote the set of all triangular
quadri-tilings of the plane, up to isometry.\jump
In \cite{Bea}, triangular quadri-tilings of $\Q$ are shown
to correspond to two superposed dimer models in the following way, see
also Figure \ref{fig2}. Define a lozenge to be a $60^{\circ}$-rhombus. 
Then triangular quadri-tilings are $2$-tilings of a family of graphs $\LL$ which are
lozenge-with-diagonals tilings of the plane, up to isometry. Indeed let $T\in\Q$ be a
triangular quadri-tiling, then on every quadri-tile of $T$ draw the edge
separating the two right triangles, this yields a lozenge-with-diagonals
tiling $L(T)$ called the {\bf underlying tiling}. Moreover, the lozenge tiling $\Ls(T)$
obtained from $L(T)$ by removing the diagonals, is a $2$-tiling of the
equilateral triangular lattice $\TT$.\\

\noindent Note that lozenge-with-diagonals tilings and the equilateral
triangular lattice are isoradial graphs. Assigning the critical weight
function to edges of their dual graphs, we deduce that triangular
quadri-tilings of $\Q$ correspond to two superposed {\em isoradial} dimer models.

\subsubsection{Height functions for triangular quadri-tilings}\label{subsubsec022}

Triangular quadri-tilings are characterized by {\em two} height
functions in the following way. Let $T\in \Q$ be a triangular quadri-tiling,
then the first height function, denoted by $h_1$, assigns
to vertices of $T$ the ``height'' of $T$ interpreted as a
$2$-tiling of its underlying lozenge-with-diagonals tiling $L(T)$. The
second height function, denoted by $h_2$, assigns to vertices of $T$
the height of $\Ls(T)$ interpreted as a $2$-tiling of $\TT$. An
example of computation is given in Section \ref{subsec23}. By means of $h_1$
and $h_2$, triangular quadri-tilings are
interpreted in \cite{Bea} as discrete random $2$-dimensional surfaces in a
$4$-dimensional space that are projected orthogonally to the plane,
i.e. in physics terminology, as random interfaces in dimension $2+2$.

\subsubsection{Statement of result}\label{subsubsec023}
 
The notion of Gibbs measure can be extended naturally to the set $\Q$ of
all triangular quadri-tilings, see Section \ref{subsec14}. In
\cite{Bea1}, we give an explicit
expression for such a Gibbs measure $\PP$, and conjecture it to be of minimal
free energy per fundamental domain among a {\em four}-parameter family
of translation invariant, ergodic Gibbs measures. Let us assume that
triangular quadri-tilings of $\Q$ are chosen with respect to the
measure $\PP$.\jump
Corollary \ref{thm2} below describes the fluctuations of the
unnormalized height functions $h_1^\eps$ and $h_2^\eps$. Suppose that
the equilateral triangular lattice has edge-lengths
$1$, and let $\TT^\eps$ be the lattice $\TT$ whose edge-lengths have
been multiplied by $\eps$. Observe that vertices of $\TT$ are
vertices of $L$, for every lozenge-with-diagonals tiling $L\in\LL$. 
For $i=1,2$, and for $\vphi\in\K$ define:
\begin{equation*}
H_i^\eps \vphi =
\eps^2 \sum_{v\in V(\TT^\eps)}\frac{\sqrt{3}}{2}\vphi(v)h_i^\eps(v),
\end{equation*}
\begin{cor}\label{thm2}
For $i=1,2$, and every 
$\vphi_1,\ldots,\vphi_k \in \K$, $(H_i^\eps
\vphi_1,\ldots,H_i^\eps \vphi_k)$ converges in law (as $\eps
\rightarrow 0$) to $\frac{1}{\sqrt{\pi}}(F_i \vphi_1,\ldots,F_i
\vphi_k)$, where $F_i$ is a Gaussian free field. Moreover, $F_1$ and
$F_2$ are independent.
\end{cor}

\subsection{Outline of the paper}\label{subsec03}

\begin{itemize}
\item Section \ref{sec1}: statement of the explicit expressions of
\cite{Bea} for the Gibbs measures $\mu$ and $\PP$, that are used in
the proof of Theorem \ref{thm1} and Corollary \ref{thm2}.
\item Section \ref{sec2}: definition of the height function on 
vertices of $2$-tilings of isoradial graphs.
\item Section \ref{sec3}: definition of the
Gaussian free field of the plane.
\item Section \ref{sec4} and Section \ref{sec5}: proof of Theorem
\ref{thm1} and Corollary \ref{thm2}.
\end{itemize}

\noindent {\em Acknowledgments:} We would like to thank Richard Kenyon for
proposing the questions solved in this paper, and for the many
enlightening discussions. We are grateful to Erwin Bolthausen, C\'edric
Boutillier and Wendelin Werner for their advice and suggestions.

\section{Minimal free energy Gibbs measure for isoradial dimer models}\label{sec1}

In the whole of this section, we let $G$ be an isoradial
graph, whose dual graph $G^*$ is bipartite; $B$ denotes the set of
black vertices, and $W$ the set of white ones. In the proof of Theorem
\ref{thm1}, we use the explicit expression of \cite{Bea} for the minimal
free energy per fundamental domain Gibbs measure $\mu$ on $2$-tilings
of $G$, and in the proof of Corollary \ref{thm2}, we use the explicit
expression of \cite{Bea} for the Gibbs measure $\PP$ on triangular
quadri-tilings. The goal of this section is to state the expressions 
for $\mu$ and $\PP$. In order to do so, we first define the critical weight
function and the Dirac operator, introduced in \cite{Kenyon3}.

\subsection{Critical weight function}\label{subsec11}

\subsubsection{Definition}\label{subsubsec111}

The following definition is taken from \cite{Kenyon3}. To each edge
$e$ of $G^*$, we associate a unit side-length rhombus $R(e)$ whose 
vertices are the vertices of $e$ and of its dual edge $e^*$
($R(e)$ may be degenerate). Let $\widetilde{R}=\cup_{e\in
G^*}R(e)$. The {\bf critical weight function} $\nu$ at the
edge $e$ is defined by $\nu(e)=2\sin\theta$, where $2\theta$ is the angle of the
rhombus $R(e)$ at the vertex it has in common with $e$; $\theta$ is
called the {\bf rhombus angle} of the edge $e$. Note that $\nu(e)$ is
the length of $e^*$.

\subsubsection{Example: critical weights for triangular quadri-tilings}\label{subsubsec112}

\begin{figure}[ht]
\begin{center}
\includegraphics[height=4cm]{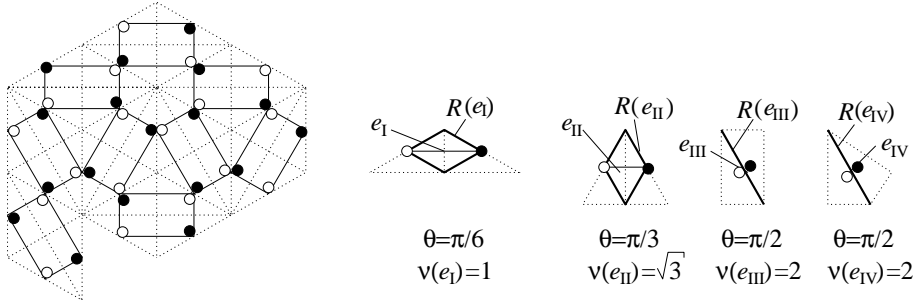}
\end{center}
\caption{Dual graph of a lozenge-with-diagonals tiling 
(left). Critical weights for quadri-tiles (right).} \label{fig3}
\end{figure}

Recall that triangular quadri-tilings of $\Q$ correspond to two superposed
isoradial dimer models, the first on lozenge-with-diagonals
tilings and the second on the equilateral triangular lattice
$\TT$. Let us note that the dual graphs of lozenge-with-diagonals tilings and
of $\TT$ are bipartite. We now compute the critical weights in the above two
cases.\jump
Consider the equilateral triangular lattice $\TT$, then
edges of its dual graph $\TT^*$, known as the honeycomb lattice, all have the same
rhombus angle, equal to $\pi/3$, and the same critical weight, equal
to $\sqrt{3}$.\jump
Consider a lozenge-with-diagonals tiling $L\in\LL$.
Observe that the circumcenters of the faces of $L$ are on the boundary of the
faces, so that in the isoradial embedding of the dual graph
$L^*$ some edges have length $0$, and the rhombi associated
to these edges are degenerate, see Figure \ref{fig3}. Since edges of $L^*$ correspond to
quadri-tiles of $L$, we classify them as being of type I, II, III,
IV. Figure \ref{fig3} below gives the rhombus angles and the critical
weights associated to edges of type I, II, III and IV, denoted by
$e_{{\rm I}},e_{{\rm II}},e_{{\rm III}},e_{{\rm IV}}$ respectively.\jump

\subsection{Dirac and inverse Dirac operator}\label{subsec12}
 
Results in this section are due to Kenyon \cite{Kenyon3}, see
also Mercat \cite{Mercat2}. Define the Hermitian matrix $K$
indexed by vertices of $G^*$ as follows. If $v_1$ and $v_2$
are not adjacent, $K(v_1,v_2)=0$. If $w \in W$ and $b \in B$ are
adjacent vertices, then $K(w,b)=\overline{K(b,w)}$ is the complex
number of modulus $\nu(w b)$ and direction pointing from $w$ to
$b$. Another useful way to say this is as follows. Let $R(wb)$ be
the rhombus associated to the edge $wb$, and denote by $w,x,b,y$
its vertices in cclw (counterclockwise) order, then $K(w,b)$ is
$i$ times the complex vector $x-y$. If $w$ and $b$ have the same
image in the plane, then $|K(w,b)|=2$, and the direction of
$K(w,b)$ is that which is perpendicular to the corresponding dual
edge, and has sign determined by the local orientation. The
infinite matrix $K$ defines the {\bf Dirac operator} $K$:
$\CC^{V(G^*)} \rightarrow \CC^{V(G^*)}$, by
\begin{equation*}\label{18}
(Kf)(v)=\sum_{u \in G^*} K(v,u)f(u),
\end{equation*}
where $V(G^*)$ denotes the set of vertices of the graph $G^*$.\jump
The {\bf inverse Dirac operator} $K^{-1}$ is defined to be
the operator which satisfies
\begin{enumerate}
\item $K K^{-1}=\mbox{Id}$, \item $K^{-1}(b,w) \rightarrow 0$,
when $|b-w| \rightarrow \infty$.
\end{enumerate}
In \cite{Kenyon3}, Kenyon proves uniqueness of $K^{-1}$, and existence
by giving an explicit expression for $K^{-1}(b,w)$ as a
function of the {\em local} geometry of the graph.

\subsection{Minimal free energy Gibbs measure for isoradial graphs}\label{subsec13}

If $e_1=w_1 b_1,\ldots,e_k=w_k b_k$ is a subset of edges of $G^*$,
define the {\bf cylinder set} $\{e_1,\ldots,e_k\}$ of $G^*$ to be the
set of dimer configurations of $G^*$ which contain
the edges $e_1,\ldots,e_k$. Let $\A$ be the field consisting of
the empty set and of the finite disjoint unions of cylinders.
Denote by $\sigma(\A)$ the $\sigma$-field generated by $\A$.
\begin{thm}{\rm\bf \cite{Bea}}\label{thm11}
Assume $G^*$ is doubly periodic. Then, there is a probability measure
$\mu$ on $(\M(G^*),\sigma(\A))$ such that for every cylinder
$\{e_1,\ldots,e_k\}$ of $G^*$,
\begin{equation}\label{41}
\mu(e_1,\ldots,e_k)=\left( \prod_{i=1}^{k} K(w_i,b_i)\right)
\det_{1 \leq i,\;j \leq k} \left(K^{-1}(b_i,w_j)\right).
\end{equation}
Moreover $\mu$ is a Gibbs measure on $\M(G^*)$, and it is the unique
Gibbs measure which has minimal free energy per fundamental domain
among the two-parameter family of translation invariant, ergodic Gibbs 
measures of \cite{Sheffield0}.
\end{thm}

\begin{rem}\label{rem1} $\,$\jump
$\bullet$ {\rm Refer to \cite{KOS} for the definition of the free
energy per fundamental domain.} \jump
$\bullet$ {\rm In \cite{Bea}, we prove that the periodicity assumption can be
released in the case of lozenge-with-diagonals tilings. That is, given
any lozenge-with-diagonals tiling of $\LL$, equation (\ref{41}) defines a Gibbs
measure on dimer configurations of its dual graph. Although
fundamental domains make no sense in case of non-periodic graphs, the minimal free
energy property can still be interpreted in some wider sense.}
\end{rem}

\subsection{Gibbs measure on triangular quadri-tilings}\label{subsec14}

The construction of this section is taken from \cite{Bea}.
Consider the set $\Q$ of all triangular quadri-tilings of the plane
up to isometry, and assume that quadri-tiles are assigned a positive weight
function. Then the notion of {\bf Gibbs measure on $\Q$} is a natural extension of
the one used in the case of dimer configurations of fixed graphs. It is a probability
measure that satisfies the following: if a triangular quadri-tiling is
fixed in an annular region, then triangular quadri-tilings inside and outside of
the annulus are independent; moreover, the probability of any interior
triangular quadri-tiling is proportional to the product of the weights
of the quadri-tiles. Denoting by $\M$ the set of dimer configurations
corresponding to triangular quadri-tilings of $\Q$, and using the
bijection between $\Q$ and $\M$, we obtain the definition of a Gibbs
measure on $\M$.\jump
Define $\LL^*$ to be set of dual graphs $L^*$ of
lozenge-with-diagonals tilings $L\in \LL$. Although some edges of
$\LL^*$ have length $0$, we think of them as edges of the one skeleton
of the graphs, so that to every edge of $\LL^*$, there corresponds a
unique quadri-tile. Let $e$ be an edge of $\LL^*$, and let $q_{e}$ be
the corresponding quadri-tile, then $q_e$ is made of two adjacent
right triangles. If the two triangles share the hypotenuse edge, they
belong to two adjacent lozenges; else if they share a leg, they belong
to the same lozenge. Let us call these lozenge(s) the {\bf lozenge(s) 
associated to the edge $e$}, and denote it/them by $\ls_e$ (that is
$\ls_e$ consists of either one or two lozenges). Let $\ks_e$ be the
edge(s) of $\TT^*$ corresponding to the lozenge(s) $\ls_e$. Let us 
introduce one more definition, if
$\{e_1,\ldots,e_k\}$ is a subset of edges of $\LL^*$, then the
{\bf cylinder set} $\{e_1,\ldots,e_k\}$ is the set of dimer configurations
of $\M$ which contain these edges. Denote by $\C$ the field
consisting of the empty set and of the finite disjoint unions of
cylinders. Denote by $\sigma(\C)$ the $\sigma$-field generated by $\C$.\jump
Consider a lozenge-with-diagonals tiling $L\in\LL$, and denote by
$\mu^L$ the minimal free energy per
fundamental domain Gibbs measure on $(\M(L^*),\sigma(\A))$ given by Theorem
\ref{thm11}, where $\sigma(\A)$ is the $\sigma$-field of cylinders of $\M(L^*)$.
Similarly, denote by $\mu^\TT$ the minimal free
energy per fundamental domain Gibbs measure on
$(\M(\TT^*),\sigma(\B))$, where $\sigma(\B)$ is the $\sigma$-field of
cylinders of $\M(\TT^*)$.\jump
Let us define $\tilde{\mu}^L$ on $(\M,\sigma(\C))$ by:
\begin{equation*}
\tilde{\mu}^L(e_1,\ldots,e_k)=\left\{
\begin{array}{ll}
\mu^L(e_1,\ldots,e_k) & \mbox{if the lozenges
$\ls_{e_1},\ldots,\ls_{e_k}$ belong to $\Ls$},\vspace{0.2cm}\\
0 & \mbox{else},
\end{array}
\right.
\end{equation*}
where we recall that $\Ls$ is the lozenge tiling obtained from the
lozenge-with-diagonals tiling $L$ by removing the diagonals.
Then, it is easy to check that $\tilde{\mu}^L$ is a probability measure
on $(\M,\sigma(\C))$. In order to simplify notations, we write $\mu^L$
for $\tilde{\mu}^L$ whenever no confusion occurs.\jump
Now, on $(\M\times\M(\TT^*),\B\times\C)$, define:
\begin{equation*}
\PP((e_1,\ldots,e_k)\times(k_1,\ldots,k_m))=\sum_{\{\Ls^*\in\M(\TT^*):k_1,\ldots,k_m\in\Ls^*\}}
\mu^{L}(e_1,\ldots,e_k)d\mu^\TT(\Ls^*).
\end{equation*}
Using Kolomogorov's extension theorem, $\PP$ extends to a probability
measure on $(\M\times\TT^*,\sigma(\B\times\C))$. Let us also denote by
$\PP$ the marginal of $\PP$ on $\M$, then in \cite{Bea}, $\PP$ is
shown to be a Gibbs measure on $\M$, and conjectured to be of minimal
free energy per fundamental domain among a {\em four}-parameter
family of translation invariant, ergodic Gibbs measures.

\section{Height functions}\label{sec2}

In the whole of this section, we let $G$ be an isoradial
graph whose dual graph $G^*$ is bipartite; as before, $B$ denotes the set of
black vertices, $W$ the set of white ones. We define the
{\em height function} $h$ on vertices of $2$-tilings of $G$, whose fluctuations are
described in Theorem \ref{thm1}. As in \cite{KOS}, see also
\cite{BL}, $h$ is defined using flows.\jump
The bipartite coloring of the vertices of $G^*$ induces an orientation of
the faces of $G$: color the dual faces of the black (white)
vertices black (white); orient the boundary edges of the black faces
cclw, the boundary edges of the white faces are
then oriented cw.

\subsection{Definition}\label{subsec21}

Let us first define a flow $\omega_0$ on the
edges of $G^*$. Consider an edge $wb$ of $G^*$, then $R(wb)$ is the
rhombus associated to $wb$, and $\theta_{wb}$ is the corresponding
rhombus angle. Define $\omega_0$ to be the white-to-black flow, which
flows by $\theta_{wb}/\pi$ along every edge $wb$ of $G^*$.

\begin{lem}\label{lem21}
The flow $\omega_0$ has divergence $1$ at every white vertex, and $-1$
at every black vertex of $G^*$.
\end{lem} 

\begin{proof}
By definition of the rhombus angle, we have\\
\begin{equation*}
\forall w\in W, \sum_{b:b\sim w}2\theta_{wb}=2\pi;\; \forall b~\in~B,
\sum_{w:w\sim b}2\theta_{wb}=2\pi.
\end{equation*}
\end{proof}
Now, consider a $2$-tiling $T$ of $G$, and let $M$ be the corresponding
perfect matching of $G^*$. Then $M$ defines a white-to-black unit flow
$\omega$ on the edges of $G^*$: flow by $1$ along every edge of $M$,
from the white vertex to the black one. The difference $\omega_0-\omega$ is a divergence free flow, which
means that the quantity of flow that enters any vertex of $G^*$ equals
the quantity of flow which exists that same vertex.\jump
We are ready for the definition of the {\bf height function} $h$. Choose a vertex
$v_0$ of $G$, and fix $h(v_0)=0$. For every other vertex $v$ of $T$,
take an edge-path $\gamma$ of $G$ from $v_0$ to $v$. If an edge $uv$
of $\gamma$ is oriented in the direction of the path, and if we denote
by $e$ its dual edge, then $h$ increases by $\omega_0(e)-\omega(e)$
along $uv$; if an edge $uv$ is oriented in the opposite direction,
then $h$ decreases by the same quantity along $uv$. As a consequence
of the fact that $\omega_0-\omega$ is a divergence free flow, the
height function $h$ is well defined.\jump
The following lemma gives a correspondence between height
functions defined on vertices of $G$, and $2$-tilings of $G$.
\begin{lem}\label{lem22}
Let $\tilde{h}$ be
an $\RR$-valued function on vertices of $G$ satisfying
\begin{enumerate}
\item[$\bullet$] $\tilde{h}(v_0)=0$,
\item[$\bullet$] $\tilde{h}(v)- \tilde{h}(u)=\omega_0(e)$ or
$\omega_0(e)-1$ for any edge $uv$ oriented from $u$ to $v$,
where $e$ denotes the dual edge of $uv$.
\end{enumerate}
Then, there is a bijection between functions $\tilde{h}$
satisfying these two conditions, and $2$-tilings of $G$.
\end{lem}
\begin{proof}
The idea of the proof closely follows \cite{EKLP}. Let $T$ be
a $2$-tiling of $G$, $M$ be the corresponding matching,
and $\omega$ be the unit white-to-black flow defined by $M$. Then,
the height function $h$ satisfies the conditions of the lemma:
consider an edge $uv$ of $G$ oriented from $u$ to $v$ and denote
by $e$ its dual edge, then $h(v)-h(u)=\omega_0(e)-\omega(e)$, and
by definition $\omega(e)=0$ or $1$.\\
Conversely, consider an $\RR$-valued function $\tilde{h}$ as in
the lemma. Let us construct a $2$-tiling $T$ whose height function
is $\tilde{h}$. Consider a black face $F$ of $G$, and let
$e_1,\ldots,e_m$ be the dual edges of its boundary edges. Then
$\sum_{i=1}^m \omega_0(e_i)=1$, so that there is exactly one
boundary edge $uv$ along which $\tilde{h}(v)-\tilde{h}(u)$ is
$\omega_0(e_i)-1$ (where $e_i$ is the dual edge of $uv$). To the
face $F$, we associate the $2$-tile of $G$ which is crossed by the
edge $uv$. Repeating this procedure for all black faces of $G$, we
obtain $T$.
\end{proof}
\vspace{-0.8cm}

\subsection{Interpretation}\label{subsec22}

The discrete interface interpretation of $2$-tilings was first given
by Thurston in the case of lozenges \cite{Thurston}. Following him, a
$2$-tiling of $G$ can be seen as a discrete $2$-dimensional surface
$S$ in a $3$-dimensional space that has been projected orthogonally to
the plane; the ``height'' of $S$ is given by the function
$h$. Stated in physics terminology, $2$-tilings of $G$ are random interfaces in dimension $2+1$. 
 
\subsection{Example: the two height functions of triangular quadri-tilings}\label{subsec23}

\begin{figure}[ht]
\begin{center}
\includegraphics[height=5cm]{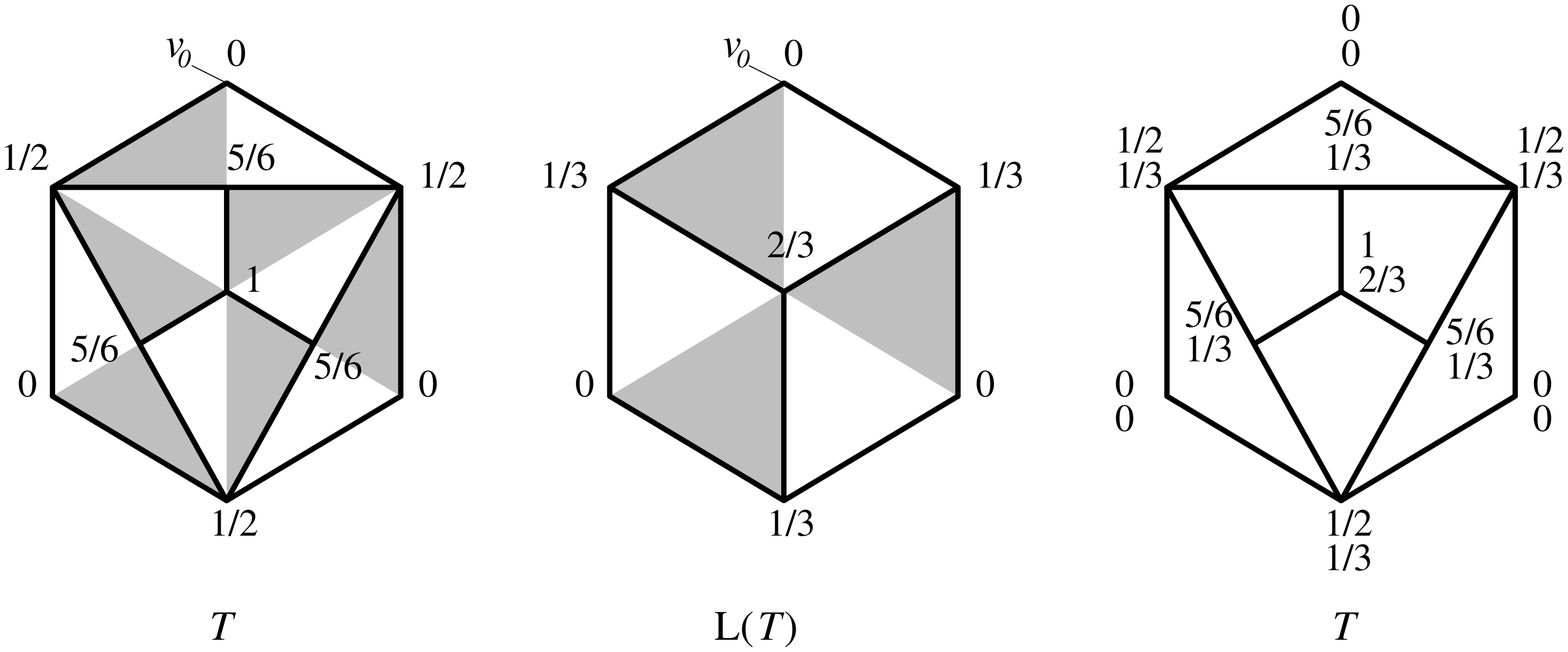}
\end{center}
\caption{From left to right: triangular quadri-tiling T and first 
height function $h_1$, underlying lozenge tiling $\Ls(T)$ and second
height function $h_2$, T with height functions $h_1$ (above) and $h_2$
(below).} \label{fig4}
\end{figure}

Consider a triangular quadri-tiling $T\in\Q$. Then, recall that $h_1$
assigns to vertices of $T$ the ``height'' of $T$ interpreted as a
$2$-tiling of its underlying lozenge-with-diagonals tiling $L(T)$, and
$h_2$ assigns to vertices of $T$ the ``height'' of $\Ls(T)$ interpreted as
a $2$-tiling of $\TT$. Let us now explicitly compute $h_1$ and
$h_2$.\jump
The definition of the flow $\omega_0(L(T))$ on edges of $L(T)^*$ uses
the rhombus angles of the edges of $L(T)^*$. These have been
computed in Section \ref{subsubsec112} and are equal to
$\frac{\pi}{6},\frac{\pi}{3},\frac{\pi}{2}$ for edges of type I, II, III and IV
respectively. Hence if $uv$ is the boundary edge of a quadri-tile 
of $T$, oriented from $u$ to $v$, then the height change of $h_1$
along $uv$ is $\frac{1}{6},\frac{1}{3},\frac{1}{2}$ depending on
whether the dual edge $e$ of the edge $uv$ is of type I, II,
III or IV respectively, see Figure \ref{fig4}.\jump
In a similar way, the flow $\omega_0(\TT^*)$ on edges of $\TT^*$ flows
by $\frac{1}{3}$ along every edge, so that if $u'v'$ be a boundary edge of a 
lozenge of $\Ls(T)$ oriented from $u'$ to $v'$, then the height
change of $h_2$ along $u'v'$ is equal to $\frac{1}{3}$. Thinking of a
lozenge as the projection to the plane of the face of a cube
\cite{Thurston}, there is a natural way to assign a value for $h_2$ at
the vertex at the crossing of the diagonals of the lozenges, see
Figure \ref{fig4}.\\

\noindent By means of $h_1$ and $h_2$, a triangular quadri-tiling $T$ 
of $\Q$ is interpreted in \cite{Bea} as a $2$-dimensional discrete surface
$S_1$ in a $4$-dimensional space that has been projected orthogonally
to the plane. $S_1$ can also be projected to
$\frac{1}{3}\widetilde{\ZZ}^3$ ($\widetilde{\ZZ}^3$ is the space
$\ZZ^3$ where the cubes are drawn with diagonals on their faces), and
one obtains a surface $S_2$. When projected to the plane $S_2$ is the
underlying lozenge-with-diagonals tiling $L(T)$. This can be restated
by saying that triangular quadri-tilings of $\Q$ are discrete interfaces in
dimension $2+2$.\newpage

\section{Gaussian free field of the plane}\label{sec3}

Theorem \ref{thm1} and Corollary \ref{thm2} prove convergence of the
height function $h$ to a Gaussian free field. The goal of this section
is to define the {\em Gaussian free field of the plane}. We refer to \cite{GJ,Sheffield1} 
for other ways of the defining the Gaussian free field.

\subsection{The Green function of the plane, and Dirichlet
energy}\label{subsec31}

\noindent The {\bf Green function of the plane}, denoted by $g$,
is the kernel of the Laplace equation in the plane, it satisfies
$\Delta_x g(x,y)=\delta_x(y)$, where $\delta_x$ is the Dirac
distribution at $x$. Up to an additive constant, $g$ is given by
\begin{equation*}\label{53}
g(x,y)=-\frac{1}{2 \pi} \log|x-y|.
\end{equation*}
Define the following bilinear form
\begin{equation*}\label{54}
\begin{array}{rccl}
G:&C_{c,0}^{\infty}(\RR^2)\times C_{c,0}^{\infty}(\RR^2) &\longrightarrow& \RR\\
&(\vphi_1,\vphi_2) &\longmapsto &
 G(\vphi_1,\vphi_2)=\displaystyle \int_{\RR^2} \int_{\RR^2} g(x,y) \vphi_1(x) \vphi_2(y)dx\,dy.\\
\end{array}
\end{equation*}
$G(\vphi,\vphi)$ is called the {\bf Dirichlet energy} of $\vphi$.
Let us consider the topology induced by the $L^{\infty}$ norm on
$\K$.

\begin{lem}\label{lem51}
$G$ is a continuous, positive definite, bilinear form.
\end{lem}
\begin{proof}
\begin{souligne}
{$G$ is continuous}
\end{souligne}
\noindent This is a consequence of the fact that for every
$\vphi_1,\vphi_2\in\K$, the function $g(x,y)$ is integrable.\jump
\begin{souligne}
{$G$ is positive definite}
\end{souligne}
\noindent For $i=1,2$, denote by $K_i=\mbox{supp}(\vphi_i)$, and
let $f_i(x)=\int_{\RR^2}g(x,y)\vphi_i(y) dy$. Let us prove that
\begin{equation}\label{55}
G(\vphi_1,\vphi_2)=\int_{\RR^2} \nabla f_1(x).\nabla f_2(x)dx.
\end{equation}
For any $R >0$, Green's formula implies
{\small\begin{equation}\label{56} \int_{B(0,R)}\nabla
f_1(x).\nabla f_2(x)dx= -\int_{B(0,R)}\triangle f_1(x)f_2(x)dx
+\int_{S(0,R)}f_2(x) \nabla f_1(x).n(x)ds.
\end{equation}}

\noindent Assume $R$ is large enough so that $K_1,K_2 \subset
B(0,R)$. The first term of the right hand side of (\ref{56})
satisfies {\small\begin{eqnarray*}\label{57}
-\int_{B(0,R)}\triangle f_1(x)f_2(x)dx&=&
-\int_{B(0,R)}\triangle_x
\left(\int_{K_1}g(x,y)\vphi_1(y)dy\right)
 \left( \int_{K_2}g(x,y)\vphi_2(y)dy \right) dx,\\
&=& \int_{K_1} \int_{K_2}g(x,y)\vphi_1(x)\vphi_2(y)dy\,dx =
G(\vphi_1,\vphi_2).
\end{eqnarray*}}

\noindent In order to evaluate the second term of the right hand
side of (\ref{56}), let us compute
\begin{eqnarray*}\label{58}
\nabla f_1(x).n(x)&=&-\frac{1}{2 \pi} \int_{K_1}\vphi_1(y)\nabla
\log|x-y|.n(x)dy,\\
&=&-\frac{1}{2\pi}\int_{K_1}\vphi_1(y)\frac{|x|}{|x-y|^2}dy,\\
&=&-\frac{1}{2\pi}\int_{K_1}\vphi_1(y)\left(
\frac{|x|}{|x-y|^2}-\frac{1}{|x|}\right) dy
-\frac{1}{2\pi}\int_{K_1}\vphi_1(y)\frac{1}{|x|}dy,\\
&=&-\frac{1}{2\pi}\int_{K_1}\vphi_1(y)\left(
\frac{|x|}{|x-y|^2}-\frac{1}{|x|}\right) dy \mbox{ (since
$\vphi_1$ is a mean $0$ function)}.
\end{eqnarray*}
$\forall \,x \in S(0,R)$, $\forall \, y \in K_1$, we have
$\frac{|x|}{|x-y|^2}-\frac{1}{|x|}=O\left(\frac{1}{R^2}\right)$,
hence $|\nabla f_1(x).n(x)|=O\left(\frac{1}{R^2}\right)$; $\forall
\,x \in S(0,R)$, we also have $|f_2(x)|=O(\log R)$, thus the
second term of the right hand side of (\ref{56}) is
$O\left(\frac{\log R}{R}\right)$. Taking the limit as $R
\rightarrow \infty$ in (\ref{56}), we
obtain (\ref{55}).\\

\noindent Let us assume $G(\vphi_1,\vphi_1)=0$. By equality
(\ref{55}) this is equivalent to $\int_{\RR^2} |\nabla
f_1(x)|^2dx=0$, hence $\nabla f_1 \equiv (0,0)$. Since
$\vphi_1(x)=\triangle f_1(x)=\mbox{div}(\nabla f_1(x))$, we deduce
$\vphi_1\equiv 0$.
\end{proof}

\subsection{Random distributions}\label{subsec32}

\noindent The following definitions are taken from \cite{GV}. A
{\bf random function} $F$ associates to every function $\vphi \in
\K$ a real random variable $F \vphi$. For $\vphi_1,\ldots,\vphi_k
\in \K$, we suppose that the joint probabilities $a_n \leq F
\vphi_n < b_n$, $1 \leq n \leq k$ are given, and we ask that they
satisfy the compatibility relation. A random function $F$ is {\bf
linear} if $\forall\, \vphi_1, \vphi_2 \in \K$,
\begin{equation*}\label{59}
F(\alpha \vphi_1 + \beta \vphi_2)= \alpha F \vphi_1 + \beta F
\vphi_2.
\end{equation*}
It is {\bf continuous} if convergence of the functions
$\vphi_{n_j}$ to $\vphi_j$ $(1 \leq j \leq k)$ implies
\begin{equation*}\label{510}
\lim_{n \rightarrow \infty} (F \vphi_{n_1},\ldots,F
\vphi_{n_k})=(F \vphi_1,\ldots,F \vphi_k),
\end{equation*}
that is, if $P(x)$ (resp. $P_n(x)$) is the probability measure
corresponding to the random variable $(F \vphi_1, \ldots,F
\vphi_k)$ (resp. $(F \vphi_{n_1},\ldots,F \vphi_{n_k})$), then for
any bounded continuous function $f$
\begin{equation*}\label{511}
\displaystyle \lim_{n \rightarrow \infty} \int f(x_1,\ldots,x_k)d
P_n(x)=\int f(x_1,\ldots,x_k)dP(x).
\end{equation*}

\noindent A {\bf random distribution} $F$ is a random function
which is linear and continuous. It is said to be {\bf Gaussian} if
for every linearly independent functions $\vphi_1,\ldots,\vphi_k
\in \K$,
the random vector $(F \vphi_1,\ldots,F \vphi_k)$ is Gaussian.\\

\noindent Two random distributions $F$ and $G$ are said to be {\bf
independent} if for any functions $\vphi_1,\ldots,\vphi_k \in \K$,
the random vectors $(F\vphi_1,\ldots,F \vphi_k)$ and $(G
\vphi_1,\ldots,G \vphi_k)$ are independent.

\subsection{Gaussian free field of the plane}\label{subsec33}

\begin{thm}{\rm\cite{Bochner}}\label{thm51}
If $G : \K \times \K \rightarrow \RR$ is a bilinear, continuous,
positive definite form, then there exists a Gaussian random
distribution $F$, whose covariance function is given by
$$
\EE(F \vphi_1 F \vphi_2)= G(\vphi_1,\vphi_2).
$$
\end{thm}

\noindent Using Lemma \ref{lem51}, and Theorem \ref{thm51}, we
define a {\bf Gaussian free field of the plane} to be a Gaussian
random distribution whose covariance function is
\begin{equation*}\label{512}
\EE(F \vphi_1 F \vphi_2)= -\frac{1}{2 \pi} \int_{\RR^2}
\int_{\RR^2}\log|x-y|\vphi_1(x) \vphi_2(y)\,dx\,dy.
\end{equation*}

\section{Proof of Theorem \ref{thm1}}\label{sec4}

We place ourselves in the context of Theorem \ref{thm1}: $G$ is an
isoradial graph whose dual graph $G^*$ is doubly periodic and bipartite;
edges of $G^*$ are assigned the critical weight function, and
$G^\eps$ is the graph $G$ whose edge-lengths have been multiplied by
$\eps$. Recall the following notations: 
$H^\eps\vphi=\eps^2\sum_{v\in V(G^\eps)}a(v^*)\vphi(v)h^\eps(v)$,
$\mu$ is the minimal free energy per fundamental domain Gibbs
measure for the dimer model on $G^*$ of Section \ref{subsec13},
and $F$ is a Gaussian free field of the plane.\jump
Since the random vector $(F \vphi_1,\ldots,F \vphi_k)$ is
Gaussian, to prove convergence of
$(H^\eps\vphi_1,\\ \ldots,H^\eps \vphi_k)$ to $(F
\vphi_1,\ldots,F \vphi_k)$, it suffices to prove convergence of
the moments of $(H^\eps\vphi_1,\ldots,$$H^\eps \vphi_k)$ to
those of $(F \vphi_1,\ldots,F \vphi_k)$; that is we need to show
that for every $k$-tuple of positive integers $(m_1, \ldots,m_k)$,
we have
\begin{equation}\label{514}
\lim_{\eps \rightarrow
0}\EE[(H^\eps\vphi_1)^{m_1}\ldots(H^\eps\vphi_k)^{m_k}]=\EE[(F
\vphi_1)^{m_1}\ldots(F \vphi_k)^{m_k}].
\end{equation}
In Section \ref{subsec41}, we prove two properties of the
height function $h$, and in Section \ref{subsec42}, we give the
asymptotic formula of \cite{Kenyon3} for the inverse Dirac operator
$K^{-1}$. Using these results in Section \ref{subsec43}, we prove 
a formula for the limit (as $\eps \rightarrow 0$) of the
$k^{\mbox{{\scriptsize th}}}$ moment of $h$. This allows us to show
convergence of
$\EE[(H^\eps\vphi)^k]$ to $\EE[(F \vphi)^k]$ in Section \ref{subsec44}. 
One then obtains equation
(\ref{514}) by choosing $\vphi$ to be a suitable linear
combination of the $\vphi_i$'s.\jump 
As before, $B$ denotes the set of black vertices of $G^*$, and $W$ the
set of white ones. Moreover, we suppose that faces of $G$ have the
orientation induced by the bipartite coloring of the vertices of $G^*$.

\subsection{Properties of the height function}\label{subsec41}

\noindent Let $u,\,v$ be two vertices of $G$, and let $\gamma$ be
an edge-path of $G$ from $u$ to $v$. First, consider edges of
$\gamma$ which are oriented in the direction of the path, that is
edges which have a black face of $G$ on the left, and denote by
$f_1,\ldots,f_n$ their dual edges. Hence an edge $f_j$ consists of
a black vertex on the left of $\gamma$, and of a white one on the
right. Similarly, consider edges of $\gamma$ which are oriented in
the opposite direction, and denote by $e_1,\ldots,e_m$ their dual
edges, hence an edge $e_i$ consists of a white vertex on the left
of $\gamma$, and of a black one on the right. Let $\II_e$ be the
indicator function of $\M(G^*)$: $\II_e(M)=1$, if the edge $e$
belongs to the dimer configuration $M$ of $G^*$, and $0$ else.

\begin{lem}\label{lem41}
\begin{equation*}\label{515}
h(v)-h(u)=\sum_{j=1}^m( \II_{e_j}-\mu(e_j))+\sum_{j=1}^n
(-\II_{f_j}+\mu(f_j)).
\end{equation*}
\end{lem}
\begin{proof}
\noindent Let $e_j$ be the dual edge of an edge $u_j v_j$ of
$\gamma$ oriented from $v_j$ to $u_j$. Denote by $\theta_j$ the
rhombus angle of the edge $e_j$, then by Lemma \ref{lem22},
\begin{equation*}\label{516}
h(v_j)-h(u_j)= \left\{
\begin{array}{ll}
-\frac{\theta_j}{\pi}+1 & \mbox{if the edge $e_j$ belongs to the dimer 
configuration of $G^*$},\vspace{0.2cm}\\
-\frac{\theta_j}{\pi} & \mbox{else}.
\end{array}
\right.
\end{equation*}
Hence
$h(v_j)-h(u_j)=(-\frac{\theta_j}{\pi}+1)\II_{e_j}-\frac{\theta_j}{\pi}(1-\II_{e_j})=\II_{e_j}-\frac{\theta_j}{\pi}$.
Moreover, a direct computation using formula (\ref{41}) yields
$\mu(e_j)=\theta_j/\pi$, so
\begin{equation*}\label{517}
h(v_j)-h(u_j)=\II_{e_j}-\mu(e_j).
\end{equation*}
Similarly, when $f_j$ is the dual edge of an edge $u_j' v_j'$ of
$\gamma$ oriented from $u_j'$ to $v_j'$, we obtain
$h(v_j')-h(u_j')=-\II_{f_j}+\mu(f_j)$, and we conclude
\begin{equation*}\label{518}
h(v)-h(u)=\sum_{j=1}^{m}h(v_j)-h(u_j)+\sum_{j=1}^n
h(v_j')-h(u_j')=\sum_{j=1}^{m}(\II_{e_j}-\mu(e_j))+\sum_{j=1}^n(-\II_{f_j}+\mu(f_j)).
\end{equation*}
\end{proof}
\begin{lem}\label{lem42}
\begin{equation*}\label{519}
\EE_{\mu}[h(v)-h(u)]=0.
\end{equation*}
\end{lem}
\begin{proof}
By Lemma \ref{lem41} we have,
$\displaystyle
\EE_{\mu}[h(v)-h(u)]=\sum_{j=1}^m
\EE_{\mu}[\II_{e_j}-{\mu}(e_j)]+\sum_{j=1}^n \EE_{\mu}
[-\II_{f_j}+{\mu}(f_j)]=0.
$
\end{proof}

\subsection{Asymptotics of the inverse Dirac operator $K^{-1}$}\label{subsec42}

In order to state the asymptotic formula of \cite{Kenyon3} for the inverse of the Dirac
operator $K$ indexed by vertices of $G^*$, we let $w\in W$ be a white vertex
of $G^*$, $v$ any other vertex of $G^*$, and define the rational
function $f_{wv}(z)$ of \cite{Kenyon3}. Recall that $\widetilde{R}$ is the set of
rhombi associated to edges of $G^*$, and consider $w=v_0,v_1,v_2,\ldots,v_k=v$ an
edge-path of $\widetilde{R}$ from $w$ to $v$. Each edge $v_j
v_{j+1}$ has exactly one vertex of $G^*$ (the other is a vertex of
$G$). Direct the edge away from this vertex if it is white, and
towards this vertex if it is black. Let $e^{i\alpha_j}$ be the
corresponding vector in $\widetilde{R}$ (which may point contrary
to the direction of the path). Then, $f_{wv}$ is defined inductively
along the path, starting from
$$f_{ww}(0)=1.$$
If the edge leads away from a white vertex, or towards a black
vertex, then
$$
f_{w v_{j+1}}(z)=\frac{f_{w v_j}(z)}{1-e^{i\alpha_j}},
$$
else, if it leads towards a white vertex, or away from a black
vertex, then
$$
f_{w v_{j+1}}(0)=f_{w v_j}(z).(z-e^{i\alpha_j}).
$$
The function $f_{wv}(z)$ is well defined (i.e. independent of the
edge-path of $\widetilde{R}$ from $w$ to $v$). Then, Kenyon gives the
following asymptotics for the inverse Dirac operator $K^{-1}$:
\begin{thm}{\rm\cite{Kenyon3}}\label{thm15}
Asymptotically, as $|b-w|\rightarrow \infty$,
$$
K^{-1}(b,w)=\frac{1}{2 \pi} \left(
\frac{1}{b-w}+\frac{f_{wb}(0)}{\bar{b}-\bar{w}}\right)+
O\left(\frac{1}{|b-w|^2}\right).
$$
\end{thm}

\subsection{Moment formula}\label{subsec43}

\noindent Let $u_1,\ldots,u_k,v_1,\ldots,v_k$ be distinct points
of $\RR^2$, and let $\gamma_1,\ldots,\gamma_k$ be pairwise
disjoint paths such that $\gamma_j$ runs from $u_j$ to $v_j$. Let
$u_j^{\eps},v_j^{\eps}$ be vertices of $\Geps$ lying within
$O(\eps)$ of $u_j$ and $v_j$ respectively. In order to simplify
notations, we write $h$ for the unnormalized height function $h^\eps$
on $2$-tilings of $G^\eps$. Then, we have

\begin{prop}\label{prop51}
For every $k \in \NN$, $k \geq 2$
{\scriptsize \begin{equation}\label{521}
\displaystyle\lim_{\eps \rightarrow 0}
\EE[(h(v_1^{\eps})-h(u_1^{\eps}))
\ldots(h(v_k^{\eps})-h(u_k^{\eps}))]=\\
\frac{(-i)^k}{(2\pi)^k}\sum_{\eps=0,1}(-1)^{k\eps}\left(
\int_{\gamma_1}\ldots\int_{\gamma_k} \det_{
\begin{array}{l}
\scriptstyle i,j \in [1,k]\\
\scriptstyle i \neq j
\end{array}} \left(
\frac{1}{z_i^\eps-z_j^\eps} \right) dz_1^\eps\ldots
dz_k^\eps\right),
\end{equation}}
\noindent where $z_i^0=z_i$ and $z_i^1=\bar{z_i}$.
\end{prop}
\begin{proof}
Steps of the proof follow \cite{Kenyon1}, but since we work in a much
more general setting, they are adapted in a non-trivial way.\jump
Let $\gamma_1^{\eps},\ldots,\gamma_k^{\eps}$ be pairwise disjoint
paths of $\Geps$, such that $\gamma_j^{\eps}$ runs from
$u_j^{\eps}$ to $v_j^{\eps}$ and approximates $\gamma_j$ within
$O(\eps)$. For every $j$, denote by $f_{js}$ the dual edge of the
$s^{\mbox{{\scriptsize th}}}$ edge of the path $\gamma_j^\eps$,
which is oriented in the direction of the path: $f_{js}$ consists
of a black vertex on the left of $\gamma_j^{\eps}$, and of a white
one on the right. Denote by $e_{jt}$ the dual edge of the
$t^{\mbox{{\scriptsize th}}}$ edge of the path $\gamma_j^\eps$,
which is oriented in the opposite direction: $e_{jt}$ consists of
a black vertex on the right of $\gamma_j^{\eps}$, and of a white
one on the left. Using Lemma \ref{lem41}, we obtain\jump
$\EE[(h(v_1^{\eps})-h(u_1^{\eps}))
\ldots(h(v_k^{\eps})-h(u_k^{\eps}))]=$
{\scriptsize\begin{eqnarray}
&=&\EE\left[\sum_{t_1}(\II_{e_{1t_1}}-{\mu}(e_{1t_1}))-\sum_{s_1}
(\II_{f_{1 s_1}}-{\mu}(f_{1 s_1}))\right]\ldots
\left[\sum_{t_k}(\II_{e_{kt_k}}-{\mu}(e_{kt_k}))-\sum_{s_k}(\II_{f_{ks_k}}-{\mu}(f_{ks_k}))\right]
,\nonumber\\
&=&\sum_{t_1,\ldots,t_k} \EE[\II_{e_{1
t_1}}-{\mu}(e_{1t_1})]\ldots[\II_{e_{k t_k}}-{\mu}(e_{k t_k})]-
\ldots +(-1)^k \sum_{s_1, \ldots,s_k} \EE[\II_{f_{1
s_1}}-{\mu}(f_{1 s_1})]\ldots
[\II_{f_{k s_k}}-{\mu}(f_{k s_k})],\nonumber\\
&=&\sum_{\delta_1,\ldots,\delta_k \in \{0,1\}}
\sum_{t_1^{\delta_1},\ldots,t_k^{\delta_k}}(-1)^{\delta_1+\ldots+\delta_k}
\EE[\II_{e_{1 t_1^{\delta_1}}}-{\mu}(e_{1
t_1^{\delta_1}})]\ldots [\II_{e_{k t_k^{\delta_k}}}-{\mu}(e_{k
t_k^{\delta_k}})],\label{522}
\end{eqnarray}}
where $t_j^{\delta_j}= \left\{
\begin{array}{l}
t_j \mbox{ if }\delta_j=0\\
s_j \mbox{ if }\delta_j=1
\end{array}
\right.,\,e_{j t_j^{\delta_j}}=e_{j t_j^{\delta_j}}^{\delta_j}=
\left\{
\begin{array}{l}
e_{j t_j} \mbox{ if }\delta_j=0\\
f_{j s_j} \mbox{ if }\delta_j=1
\end{array}
\right.. $\jump 
For the time being, let us drop the second subscript. Write
$e_j=w_j b_j$ and $f_j=w_j' b_j'$. Moreover, let us
introduce the notation $w_j^{\delta_j}$, where
$w_j^{\delta_j}=w_j$ if $\delta_j=0$, and $w_j^{\delta_j}=w_j'$ if
$\delta_j=1$, similarly we introduce the notation
$b_j^{\delta_j}$. Hence we can write a generic term of (\ref{522})
as
\begin{equation*}
(-1)^{\delta_1+\ldots+\delta_k} \EE[(\II_{w_1^{\delta_1}
b_1^{\delta_1}}-{\mu}(w_1^{\delta_1}
b_1^{\delta_1}))\ldots(\II_{w_k^{\delta_k}b_k^{\delta_k}}-{\mu}(w_k^{\delta_k}b_k^{\delta_k}))].
\end{equation*}

\begin{lem} {\rm\cite{Kenyon0}}\label{lem54}
\begin{equation}\label{523}
\begin{split}
&(-1)^{\delta_1+\ldots+\delta_k} \EE[(\II_{w_1^{\delta_1}
b_1^{\delta_1}}-{\mu}(w_1^{\delta_1}
b_1^{\delta_1}))\ldots(\II_{w_k^{\delta_k}b_k^{\delta_k}}-{\mu}(w_k^{\delta_k}b_k^{\delta_k}))]=\\
&= (-1)^{\delta_1+\ldots+\delta_k}a_E \left|
\begin{array}{cccc}
0 & K^{-1}(b_1^{\delta_1},w_2^{\delta_2})& \dots & K^{-1}(b_1^{\delta_1},w_k^{\delta_k})\\
K^{-1}(b_2^{\delta_2},w_1^{\delta_1})&0&&\vdots\\
\vdots&&&K^{-1}(b_{k-1}^{\delta_{k-1}},w_k^{\delta_k})\\
K^{-1}(b_k^{\delta_k},w_1^{\delta_1})&\dots&K^{-1}(b_k^{\delta_k},w_{k-1}^{\delta_{k-1}})&0
\end{array}
\right|,
\end{split}
\end{equation}
where $\displaystyle a_E=\prod_{j=1}^{k}K(w_j^{\delta_j},
b_j^{\delta_j})$, and $K$ is the Dirac operator indexed by
 vertices of $G^*$.
\end{lem}

\noindent A typical term in the expansion of (\ref{523}) is
\begin{equation}\label{524}
(-1)^{\delta_1+\ldots+\delta_k} a_E \mbox{ {\rm sgn}}\,\sigma
K^{-1}(b_1^{\delta_1}, w_{\sigma(1)}^{\delta_{\sigma(1)}})\ldots
K^{-1}(b_k^{\delta_k},w_{\sigma(k)}^{\delta_{\sigma(k)}}),
\end{equation}
where $\sigma \in \widetilde{S_k}$, and $\widetilde{S_k}$ is the
set of permutations of $k$ elements, with no fixed points. To
simplify notations, let us assume $\sigma$ is a $k$-cycle, hence
(\ref{524}) becomes
\begin{equation}\label{525}
(-1)^{\delta_1+\ldots+\delta_k} a_E \mbox{ {\rm sgn}}\,\sigma
K^{-1}(b_1^{\delta_1}, w_{2}^{\delta_2})\ldots
K^{-1}(b_k^{\delta_k},w_{1}^{\delta_1}).
\end{equation}

\begin{lem}\label{lem55}
When $\eps$ is small, and for every $\delta_1,\ldots,\delta_k \in
\{0,1\}$,
\begin{equation}\label{526}
\eps^k a_E =(-i)^k
(-1)^{\delta_1+\ldots+\delta_k}dz_1^{\delta_1}\ldots
dz_k^{\delta_k}.
\end{equation}
\end{lem}
\begin{proof}
Let $u_j v_j$ be an edge of the path $\gamma_j^\eps$ where $u_j$
precedes $v_j$. We can write
\begin{equation}\label{527}
u_j v_j=\eps \l(u_j v_j)e^{i \theta_j},
\end{equation}
where $\l(u_j v_j)$ is the length of the edge $u_j v_j$ in $G$,
and $\theta_j$ is the direction from $u_j$ to $v_j$. Let us first
consider the case of an edge $u_j v_j$ oriented in the direction
of the path, that is the dual edge $w_j' b_j'$ of $u_j v_j$ has
its black vertex on the left of $\gamma_j^\eps$. By definition of
the Dirac operator, we have $K(w_j',b_j')=\l(u_j v_j)e^{i
\theta_j}e^{i \frac{\pi}{2}}$. Next we consider the case of an
edge $u_j v_j$ oriented in the opposite direction, that is its
dual edge $w_j b_j$ has its black vertex on the right of
$\gamma_j^\eps$. Again, using the definition of the Dirac
operator, we obtain $K(w_j,b_j)=\l(u_j v_j)e^{i \theta_j}e^{-i
\frac{\pi}{2}}$. We summarize the two cases by the following
equation
\begin{equation}\label{528}
K(w_j^{\delta_j},b_j^{\delta_j})=(-i)(-1)^{\delta_j}\l(u_j
v_j)e^{i\theta_j}.
\end{equation}
When $\eps$ is small we replace $u_j v_j$ by $dz_j^{\delta_j}$.
Thus combining equations (\ref{527}) and (\ref{528}) we obtain
equation (\ref{526}).
\end{proof}
\begin{lem}\label{lem56}
When $\eps$ is small and up to a term of order $O(\eps)$, equation
$(\ref{525})$ equals
{\scriptsize
\begin{equation}\label{529}
\frac{(-i)^k}{(2\pi)^k} \sum_{\eps_1,\dots,\eps_k \in \{0,1\}}
 [f_{w_2^{\delta_2} b_1^{\delta_1}}(0)]^{\eps_1}\ldots
 [f_{w_1^{\delta_1} b_k^{\delta_k}}(0)]^{\eps_k}
F_{\eps_1}(b_1^{\delta_1},w_2^{\delta_2})\ldots
F_{\eps_k}(b_k^{\delta_k},w_1^{\delta_1})dz_1^{\delta^1}\ldots
dz_k^{\delta^k},
\end{equation}}
\vspace{-0.4cm}

\noindent where $\displaystyle
F_0(z,w)=\frac{1}{z-w},\,F_1(z,w)={F_0(\bar{z},\bar{w})}$, and the
functions $f_{wb}$ are defined in Section $\ref{subsec42}$.
\end{lem}
\begin{proof}
Let us drop the superscripts $\delta_i$. Plugging relation
(\ref{526}) in (\ref{525}), we obtain
\begin{equation}\label{530}
(\ref{525})=(-i)^k \mbox{sgn}\,\sigma K^{-1}(b_1,w_2)\ldots
K^{-1}(b_k,w_1)\frac{1}{\eps^k}dz_1\ldots dz_k.
\end{equation}
Moreover, for every $i \neq j$, $\displaystyle \lim_{\eps
\rightarrow 0}\frac{|b_i-w_j|}{\eps}=\infty$, so that by Theorem
\ref{thm15} we have
\begin{equation}\label{531}
K^{-1}(b_i,w_j)=\frac{\eps}{2 \pi} (F_0(b_i,w_j)+ f_{w_j
b_i}(0)F_1(b_i,w_j))+O(\eps^2).
\end{equation}
Equation (\ref{529}) is then (\ref{530}) where the elements
$K^{-1}(b_i,w_j)$ have been replaced by (\ref{531}) and expanded
out.
\end{proof}
In what follows, all that we say is true whether the edge
$w_j^{\delta_j}b_j^{\delta_j}$ has its black vertex on the right
or on the left of the path $\gamma_j^\eps$, that is whether
$\delta_j=0$ or $1$. So to simplify notations, let us write
$\{t_j\}$ instead of $\{\delta_j \in \{0,1\},t_j^{\delta_j}\}$,
hence $\{t_j\}$ is the set of indices which run along the path
$\gamma_j^\eps$. Keeping in mind that our aim is to take the limit
as $\eps \rightarrow 0$, we replace the vertices $b_j$ and $w_j$
in the argument of the function $F_{\eps_j}$ by one common vertex
denoted by $z_j$. Define\jump $H(\eps_1,\ldots,\eps_k)$=
{\scriptsize
\begin{equation*}\label{532}
=\sum_{t_1,\ldots,t_k}
 [f_{w_{2 t_2} b_{1 t_1}}(0)]^{\eps_1} \ldots
 [f_{w_{1 t_1} b_{k t_k}}(0)]^{\eps_k}
 F_{\eps_1}(z_{1 t_1},z_{2 t_2})\ldots F_{\eps_k}(z_{k t_k},z_{1 t_1})
 dz_{1 t_1}\ldots dz_{k t_k}.
\end{equation*}}

\begin{lem}\label{lem57}\
\begin{description}
\item[{\rm 1.}] If $(\eps_1,\ldots,\eps_k)=(0,\ldots,0)$, then
$$
\lim_{\eps \rightarrow 0} H(0,\ldots,0)=
\int_{\gamma_1}\ldots \int_{\gamma_k}F_0(z_1,z_2)\ldots F_0(z_k,z_1)dz_1\ldots dz_k.
$$
\item[{\rm 2.}] If $(\eps_1,\ldots,\eps_k)=(1,\ldots,1)$, then
$$
\lim_{\eps \rightarrow 0} H(1,\ldots,1)=(-1)^k
\int_{\gamma_1}\ldots \int_{\gamma_k}
F_0(\bar{z}_1,\bar{z}_2)\ldots F_0(\bar{z}_k,\bar{z}_1)
d\bar{z}_1\ldots d\bar{z}_k.
$$
\item[{\rm 3.}] Assume there exists $i \neq j \in \{1,\ldots,k\}$
such that $\eps_i=0,\,\eps_j=1$, then
$$\lim_{\eps \rightarrow 0} |H(\eps_1,\ldots,\eps_k)|=0.$$
\end{description}
\end{lem}
\begin{proof}
Here are some preliminary notations. Dropping the second
subscript, we consider an edge $u_j v_j$ of one of the paths
$\gamma_i^\eps$, where $u_j$ precedes $v_j$. Let us denote by $w_j
b_j$ the dual edge of the edge $u_j v_j$, and let $\eps
e^{i\alpha_j}=v_j-w_j$, $\eps e^{i\beta_j}=b_j-v_j$. With these
notations, we have
$dz_j=\eps(e^{i\alpha_j}-e^{i\beta_j})$ (see Figure \ref{fig51}).\\

\begin{figure}[ht]
\begin{center}
\includegraphics[height=2.6cm]{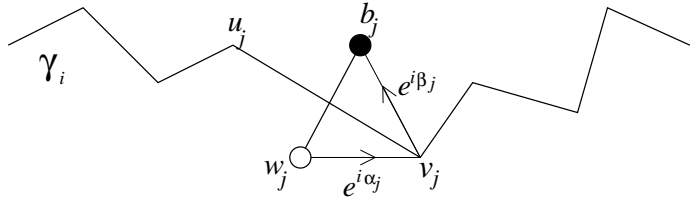}
\end{center}
\caption{Notations.} \label{fig51}
\end{figure}

\noindent Moreover, define
\begin{equation*}
J(\eps_1,\ldots,\eps_k)=
 [f_{w_2 b_1}(0)]^{\eps_1} \ldots
 [f_{w_1 b_k}(0)]^{\eps_k} dz_1\ldots dz_k.
\end{equation*}

\begin{souligne}
{\sf Proof of $1.$}
\end{souligne}
$$
J(0,\ldots,0)=dz_1\ldots dz_k,
$$
so that
$$
H(0,\ldots,0)=\sum_{t_1,\ldots,t_k}F_0(z_{1 t_1},z_{2 t_2})\ldots
F_0(z_{k t_k},z_{1 t_1})dz_{1 t_1}\ldots dz_{k t_k}.
$$
Since the paths $\gamma_j$ are disjoint, the function
$F_0(z_1,z_2)\ldots F_0(z_k,z_1)$ is integrable, and taking the
limit as $\eps \rightarrow 0$, we obtain $1$.\jump

\begin{souligne}
{\sf Proof of $2.$}
\end{souligne}
\begin{equation*}
J(1,\ldots,1)=
 f_{w_2 b_1}(0)\ldots f_{w_1 b_k}(0)\,dz_1\ldots dz_k.
\end{equation*}
Fix a vertex $v$ of $G^*$. Then, by definition of the function
$f_{wv}$,
\begin{eqnarray*}
J(1,\ldots,1)&=&
 f_{w_2 v}(0) f_{v b_1}(0)\ldots f_{w_1 v}(0) f_{v
 b_k}(0)\,dz_1\ldots dz_k,\\
&=& f_{w_1 v}(0) f_{v b_1}(0)\ldots f_{w_k v}(0) f_{v
 b_k}(0)\,dz_1 \ldots dz_k,\\
&=& f_{w_1 b_1}(0)\ldots f_{w_k b_k}(0)\, dz_1\ldots dz_k.
\end{eqnarray*}
For every $j$, we have $f_{w_j b_j}(0)=e^{-i(\beta_j+\alpha_j)}$.
Moreover, recall that $dz_j=\eps (e^{i \alpha_j}-e^{i\beta_j})$,
so that $-d\bar{z}_j=e^{-i(\beta_j+\alpha_j)}dz_j$, and we deduce
$$
J(1,\ldots,1)=(-1)^kd\bar{z}_1\ldots d\bar{z}_k.
$$
This implies,
$$
H(1,\ldots,1)=(-1)^k\sum_{t_1,\ldots,t_k} F_0(\bar{z}_{1
t_1},\bar{z}_{2 t_2})\ldots F_0(\bar{z}_{k t_k},\bar{z}_{1 t_1})
d\bar{z}_{1 t_1}\ldots d\bar{z}_{k t_k}.
$$
Taking the limit as $\eps \rightarrow 0$, we obtain $2$.\\

\begin{souligne}
{\sf Proof of $3.$}
\end{souligne}
\noindent Consider $0<\l<k$, and assume
$\eps_1=\ldots=\eps_{\l-1}=0$, $\eps_\l=\ldots=\eps_k=1$. Let us
prove that $\lim_{\eps \rightarrow 0}
|H(0,\ldots,0,1,\ldots,1)|=0.$ Note that up to a permutation of
indices, the argument is the same for the other cases.
\begin{equation*}
J(0,\ldots,0,1,\ldots,1)=
 f_{w_{\l+1} b_\l}(0) \ldots f_{w_1 b_k}(0)\,dz_1\ldots dz_k.
\end{equation*}
As above, let $v$ be a vertex of $G^*$. Then,
\begin{equation*}
J(0,\ldots,0,1,\ldots,1)=
 (f_{v b_\l}(0) dz_\l)(f_{w_1 v}(0)dz_1)
 dz_2\ldots dz_{\l-1}d\bar{z}_{\l+1}\ldots d\bar{z}_k.
\end{equation*}
Introducing the following notation
\begin{eqnarray*}
H_1&=&F(z_{2 t_2},z_{3 t_3})\ldots F(z_{{\l-2} t_{\l-2}},z_{{\l-1} t_{\l-1}})
F(\bar{z}_{{\l+1} t_{\l+1}},\bar{z}_{{\l+2} t_{\l+2}})\ldots
F(\bar{z}_{{k-1} t_{k-1}},\bar{z}_{k t_{k}}),\\
H_2&=&F(z_{1 t_1},z_{2 t_2})F(\bar{z}_{k t_k},\bar{z}_{1 t_{1}}),\\
H_3&=&F(z_{{\l-1}t_{\l-1}},z_{\l t_\l})F(\bar{z}_{\l
t_\l},\bar{z}_{\l+1 t_{\l+1}}),
\end{eqnarray*}
we obtain $H(0,\ldots,0,1,\ldots,1)=$ {\scriptsize
\begin{equation*}
=\sum_{t_2,\ldots,\hat{t_\l},\ldots,t_k}H_1
 \left(\sum_{t_1}H_2 f_{w_{1 t_1}v}(0) dz_{1 t_1}\right)
 \left(\sum_{t_\l} H_3 f_{v b_{\l t_\l}}(0)dz_{\l t_\l}\right)
 dz_{2 t_2}\ldots dz_{\l-1 t_{\l-1}}d\bar{z}_{\l+1 t_{\l+1}}\ldots d\bar{z}_{k t_k}.\\
\end{equation*}}
Let us prove
\begin{equation}\label{51}
\displaystyle \sum_{t_\l}f_{v b_{\l t_\l}}(0)dz_{\l t_\l}=O(\eps).
\end{equation}
Dropping the second subscript, let $u_1,
v_1=u_2,v_2=u_3,\ldots,v_{m-1}=u_m,v_m$ be the edge-path
$\gamma_\l^\eps$. Denote by $\xi$ the quantity $f_{v u_1}(0)$,
then
$$
f_{v b_j}(0)=\xi e^{i(\beta_1-\alpha_1)}\ldots
e^{i(\beta_{j-1}-\alpha_{j-1})} e^{-i\alpha_j}.
$$
Since $dz_j=\eps(e^{i \alpha_j}-e^{i \beta_j})$, we obtain
$$
f_{v b_j}(0)dz_j=\eps\xi\left[\left(e^{i(\beta_1-\alpha_1)} \ldots
e^{i(\beta_{j-1}-\alpha_{j-1})}\right)-\left(e^{i(\beta_1-\alpha_1)}
\ldots
e^{i(\beta_{j-1}-\alpha_{j-1})}e^{i(\beta_j-\alpha_j)}\right)\right],
$$
hence { \scriptsize
\begin{eqnarray*}
\sum_{t_\l}f_{v b_{\l t_\l}}(0)dz_{\l t_\l}&=&
 \sum_{j=1}^m f_{v b_j}(0)dz_j,\\
 &=&\eps\xi [1-e^{i(\beta_1-\alpha_1)}+\sum_{j=2}^m e^{i(\beta_1-\alpha_1)}\ldots
 e^{i(\beta_{j-1}-\alpha_{j-1})})-(e^{i(\beta_1-\alpha_1)}\ldots
 e^{i(\beta_{j-1}-\alpha_{j-1})}e^{i(\beta_j-\alpha_j})],\\
 &=&\eps \xi [1-(e^{i(\beta_1-\alpha_1)}\ldots e^{i(\beta_m-\alpha_m)}], \mbox{ (telescopic sum).}
\end{eqnarray*}}
\vspace{-0.4cm}

\noindent We deduce $|\sum_{t_\l}f_{v b_{\l t_\l}}(0)dz_{\l t_\l}|
\leq 2 \eps$,
and (\ref{51}) is proved.\\
In a similar way we prove
$\sum_{t_1}f_{w_{1 t_1}v}(0)dz_{1 t_1}=O(\eps)$.\\
Using Taylor expansion in $\eps$ for $H_2$ and $H_3$, we deduce
that $\left(\sum_{t_1}H_2 f_{w_{1 t_1}v}(0) dz_{1 t_1}\right)$ and
$\left(\sum_{t_\l} H_3 f_{v b_{\l t_\l}}(0)dz_{\l t_\l}\right)$
are $O(\eps)$. Since the function $H_1$ is integrable, we conclude
that $H(0,\ldots,0,1,\ldots,1)$ is $O(\eps^2)$ and so $3.$ is
proved.
\end{proof}
Rewriting the second subscript, and summing equation (\ref{525})
over the paths $\gamma_1^\eps,\ldots,\gamma_k^\eps$, we obtain (by
Lemmas \ref{lem56} and \ref{lem57}):
\begin{equation}\label{533}
\begin{split}
\lim_{\eps \rightarrow 0}\,\mbox{{\rm sgn}}\,\sigma&
\sum_{\delta_1,\ldots,\delta_k \in \{0,1\}}\sum_{t_1^{\delta_1}
\ldots t_k^{\delta_k}} (-1)^{\delta_1+\ldots+\delta_k} a_E
K^{-1}(b_{1 t_1^{\delta_1}},w_{2 t_2^{\delta_2}})\ldots
K^{-1}(b_{k t_k^{\delta_k}},w_{1 t_1^{\delta_1}})=\\
&=\frac{(-i)^k}{(2 \pi)^k} \mbox{{\rm
sgn}}\,\sigma\sum_{\eps=0,1}(-1)^{\eps k} \left(
\int_{\gamma_1}\ldots\int_{\gamma_k} F_0(z_1^\eps,z_2^\eps)\ldots
F_0(z_k^\eps,z_1^\eps)dz_1^\eps \ldots dz_k^\eps\right),
\end{split}
\end{equation}
where $z_i^0=z_i$ and $z_i^1=\bar{z_i}$. When $\sigma$ is a
product of disjoint cycles, we can treat each cycle separately and
the result is the product of terms like (\ref{533}). Thus when we
sum over all permutations with no fixed points, we obtain equation
(\ref{521}) of Proposition \ref{prop51}.
\end{proof}
\begin{prop}\label{prop52}\
\begin{description}
\item[-] When $k$ is odd, $\displaystyle\lim_{\eps \rightarrow 0}
\EE[(h(v_1^\eps)-h(u_1^\eps))\ldots(h(v_k^\eps)-h(u_k^\eps))]=0.$
\item[-] When $k$ is even, $\displaystyle\lim_{\eps \rightarrow 0}
\,\EE[(h(v_1^\eps)-h(u_1^\eps))\ldots(h(v_k^\eps)-h(u_k^\eps))]=$
\begin{equation*}
=\left(\frac{1}{\pi}\right)^{k/2} \sum_{\tau \in \T_k}
g(u_{\tau(1)},v_{\tau(1)},u_{\tau(2)},v_{\tau(2)}) \ldots
g(u_{\tau(k-1)},v_{\tau(k-1)},u_{\tau(k)},v_{\tau(k)}),
\end{equation*}
where $g(u,v,u',v')=g(v,v')+g(u,u')-g(v,u')-g(u,v')$, $g$ is the
Green function of the plane, and $\T_k$ is the set of all
$(k-1)!!$ pairings of $\{1,\ldots,k\}$.
\end{description}
\end{prop}
\begin{proof}
Let us cite the following lemma from \cite{Kenyon1}.
\begin{lem}{\rm\cite{Kenyon1}}\label{lem58}
Let $M=(m_{ij})$ be the $k \times k$ matrix defined by $m_{ii}=0$,
and $m_{ij}=\frac{1}{x_i - x_j}$, when $i \neq j$. Then when $k$
is odd, $\det M=0$, and when $k$ is even
$$
\det M = \sum_{\tau \in \T_k} \frac{1}{(x_{\tau(1)}-x_{\tau(2)})^2\ldots(x_{\tau(k-1)}-x_{\tau(k)})^2}.
$$
\end{lem}
Combining Proposition \ref{prop51} and Lemma \ref{lem58}, when
$k=2$, we obtain
\begin{equation*}
\begin{split}
\lim_{\eps \rightarrow 0} \EE[(h(v_1^\eps)-h(u_1^\eps)&)(h(v_2^\eps)-h(u_2^\eps))]=\\
&=-\frac{1}{4 \pi^2}\left(\int_{\gamma_1}\int_{\gamma_2}
\frac{1}{(z_1-z_2)^2}dz_1 dz_2+
\int_{\gamma_1}\int_{\gamma_2} \frac{1}{(\bar{z}_1-\bar{z}_2)^2}d\bar{z}_1 d\bar{z}_2\right),\\
&= -\frac{1}{2 \pi^2} \log \left| \frac{(v_1-v_2)(u_1-u_2)}{(v_1-u_2)(u_1-v_2)}\right|,\\
&= \frac{1}{\pi} g(u_1,v_1,u_2,v_2).
\end{split}
\end{equation*}
The case of a general even $k$ is an easy but notationally
cumbersome extension of the case $k=2$.
\end{proof}

\subsection{Proof of Theorem \ref{thm1}}\label{subsec44}

\begin{prop}\label{prop53}\
\begin{equation}
\lim_{\eps \rightarrow 0} \EE[(H^\eps
\vphi)^k]=\frac{1}{\pi^k}\EE[(F \vphi)^k]=\left\{
\begin{array}{ll}
\displaystyle 0 &\mbox{when $k$ is odd},\\
\displaystyle (k-1)!!\frac{1}{\pi^{k/2}}G(\vphi,\vphi)^{k/2}
&\mbox{when $k$ is even}.
\end{array}
\right.
\end{equation}
\end{prop}
\begin{proof}
The second equality is just the $k^{\mbox{{\scriptsize th}}}$
moment of a mean $0$, variance
$\frac{1}{\sqrt{\pi}}G(\vphi,\vphi)$, Gaussian variable. So let us
prove equality between the first and the last term.\\
Consider $u_1,\ldots, u_k$ distinct points of $\RR^2$, and for
every $j$, let $u_j^\eps$ be a vertex of $\Geps$ lying within
$O(\eps)$ of $u_j$. Define
$$
H_{u_j}^\eps \vphi=\sum_{v^\eps \in \Geps} \eps^2 a(v^*)
\vphi(v^\eps)(h(v^\eps)-h(u_j^\eps))=\sum_{v^\eps \in K^\eps}
\eps^2 a(v^*) \vphi(v^\eps)(h(v^\eps)-h(u_j^\eps)),
$$
where $K^\eps=\Geps\cap K$, and $K=$ supp$(\vphi)$, then since we
sum over a finite number of vertices,
{\scriptsize
\begin{eqnarray}
&& \EE[H_{u_1}^\eps \vphi \ldots H_{u_k}^\eps \vphi]
 = \EE \left[\sum_{v_1^\eps \in K^\eps}\eps^2
a(v_1^*)\vphi(v_1^\eps)(h(v_1^\eps)-h(u_1^\eps)) \ldots
\sum_{v_k^\eps \in K^\eps}\eps^2 a(v_k^*)\vphi(v_k^\eps)(h(v_k^\eps)-h(u_k^\eps))\right],\nonumber\\
&& = \sum_{v_1^\eps \in K^\eps} \ldots\sum_{v_k^\eps \in
K^\eps}(\eps^2)^k a(v_1^*)\ldots a(v_k^*) \vphi(v_1^\eps)\ldots \vphi(v_k^\eps)\,
\EE[(h(v_1^\eps)-h(u_1^\eps))\ldots
(h(v_k^\eps)-h(u_k^\eps))].\label{52}
\end{eqnarray}}

\begin{lem}\label{lem59}
As $\eps \rightarrow 0$, the Riemann sum {\rm (\ref{52})}
converges to
$$
\int_{\RR^2}\ldots \int_{\RR^2}\vphi(v_1)\ldots\vphi(v_k)\lim_{\eps \rightarrow 0}
\EE[(h(v_1^\eps)-h(u_1^\eps))\ldots
(h(v_k^\eps)-h(u_k^\eps))]dv_1 \ldots dv_k,
$$
where $\lim_{\eps \rightarrow 0}
\EE[(h(v_1^\eps)-h(u_1^\eps))\ldots
(h(v_k^\eps)-h(u_k^\eps))]$ is given by Proposition
$\ref{prop52}$.
\end{lem}
\begin{proof}
In what follows, all that we say is true whether $\delta_j=0$ or
$1$, so to simplify notations, as before, let us write $\{t_j\}$
instead of $\{\delta_j \in \{0,1\}, t_j^{\delta_j}\}$, hence
$\{t_j\}$ is the set of indices which run along the path
$\gamma_j^\eps$. Combining equations (\ref{522}), (\ref{524}) and
(\ref{526}) yields, \jump 
$\EE[h(v_1^\eps)-h(u_1^\eps)]\ldots
[h(v_k^\eps)-h(u_k^\eps)]=${\small
\begin{equation}\label{542}
=(-i)^k\sum_{t_1,\ldots,t_k}\sum_{\sigma \in \widetilde{S_k}}
\mbox{ sgn}\,\sigma
 K^{-1}(b_{1 t_1},w_{\sigma(1) t_{\sigma(1)}})\ldots
 K^{-1}(b_{k t_k},w_{\sigma(k) t_{\sigma(k)}})\frac{1}{\eps^k}
dz_{1 t_1}\ldots dz_{k t_k}.
\end{equation}}

\noindent It suffices to consider the case where $\sigma$ is a
$k$-cycle, other cases are treated similarly. Indices are denoted
cyclically (i.e. $k+1\equiv 1$). There is a singularity in
(\ref{542}) as soon as $v_j^\eps=v_{j+1}^\eps$ for some indices
$j$. Hence, we need to prove that for $\eps$ small enough,
{\scriptsize
\begin{equation}\label{534}
\sum_{v_1^\eps}\ldots\sum_{v_k^\eps}(\eps^2)^k a(v_1^\eps)\ldots a(v_k^\eps)
 |\vphi(v_1^\eps)\ldots \vphi(v_k^\eps)|
 \sum_{t_1,\ldots,t_k}|K^{-1}(b_{1 t_1},w_{2 t_2})\ldots K^{-1}(b_{k t_k},w_{1 t_1})|
 \frac{1}{\eps^k} dz_{1 t_1}\ldots dz_{k t_k},
\end{equation}}

\noindent is $o(1)$, when the sum is over vertices
$v_1^\eps,\ldots,v_k^\eps$ that satisfy
\begin{equation*}
|v_1^\eps-v_2^\eps|\leq\delta,\ldots,|v_m^\eps-v_{m+1}^\eps|\leq\delta,
|v_{m+1}^\eps-v_{m+2}^\eps|>\delta,\ldots,|v_k^\eps-v_1^\eps|>\delta,
\end{equation*}
for some $1\leq m < k-1$. Let $0<\beta<1$, up to a renaming of
indices, this amounts to considering vertices
$v_1^\eps,\ldots,v_k^\eps$ in $\Theta_1\cap\Theta_2\cap\Theta_3$,
where
\begin{eqnarray*}
\Theta_1&=&\{v_1^\eps,\ldots,v_m^\eps\,|\,
 \{1\leq i\leq m,\,0\leq|v_i^\eps -v_{i+1}^\eps|
 \leq\eps^\beta\},\\
\Theta_2&=&\{v_{m+1}^\eps,\ldots,v_n^\eps\,|\,
 m+1\leq j\leq n,\,\eps^\beta \leq|v_j^\eps -v_{j+1}^\eps|
 \leq\delta\},\\
\Theta_3&=&\{v_{n+1}^\eps,\ldots,v_k^\eps\,|\,
 n+1\leq \l \leq k,\, |v_\l^\eps-v_{\l+1}^\eps|>
 \delta\},
\end{eqnarray*}
for some $1\leq m<n<k-1$.\\
Since $u_1,\ldots,u_k$ are distinct vertices of $\RR^2$, and since
equation (\ref{534}) does not depend on the path $\gamma_j^\eps$
from $u_j^\eps$ to $v_j^\eps$, let us choose the paths
$\gamma_1^\eps,\ldots,\gamma_k^\eps$ as follows. Note that it
suffices to consider the part of the path $\gamma_j^\eps$ where
the vertices $b_{j t_j}$ and $w_{j t_j}$ are within distance
$\delta$ from $v_j^\eps$. Let us write $|t_j|\leq\delta$ to denote
indices $t_j$ which refer to vertices of $\gamma_j^\eps$ that are
at distance at most $\delta$ from $v_j^\eps$. Take $\gamma_j^\eps$
to approximate a straight line within $O(\eps)$, from $v_j^\eps$
to $v_j^\eps+\delta$. Moreover, ask that if one continues the
lines of $\gamma_j^\eps$ and $\gamma_{j+1}^\eps$ away from
$u_j^\eps$ and $u_{j+1}^\eps$, they intersect and form an angle
$\theta_j$. Let us use the definition of the paths
$\gamma_j^\eps$, and consider the three following cases. Whenever
it is not confusing, we shall drop the second subscript. $C$
denotes a generic constant, $A\sim B$ means $A$ and $B$ are of the
same order.
 \jump $\bullet$ $v_1^\eps,\ldots,v_m^\eps\in\Theta_1$.
 \jump By Remark \ref{rem510} below, we have $|K^{-1}(b_i,w_{i+1})|\leq
C$. This implies,
\begin{equation*}
\sum_{|t_1|\leq\delta,\ldots,|t_m|\leq\delta}
 |K^{-1}(b_{1 t_1},w_{2 t_2})\ldots K^{-1}(b_{m t_m},w_{m+1 t_{m+1}})|
 \frac{1}{\eps^m} dz_{1 t_1}\ldots dz_{m t_m}\leq
 \left(\frac{C\delta}{\eps}\right)^m.
\end{equation*}
$\bullet$ $v_{n+1}^\eps,\ldots,v_k^\eps\in\Theta_3$.
 \jump By definition of the paths $\gamma_\l^\eps$, $|b_\l-w_{\l+1}|>\delta$, so that $\lim_{\eps\rightarrow
0}\frac{|b_\l-w_{\l+1}|}{\eps}=\infty$. Using Theorem \ref{thm15}
yields
$\frac{1}{\eps}|K^{-1}(b_\l,w_{\l+1})|=O\left(\frac{1}{|b_{\l}-w_{\l+1}|}\right)\leq
\frac{C}{\delta}.$ Hence,{\scriptsize
\begin{equation*}
\sum_{|t_{n+1}|\leq\delta,\ldots,|t_k|\leq\delta}
 |K^{-1}(b_{n+1 t_{n+1}},w_{n+2 t_{n+2}})\ldots K^{-1}(b_{k t_k},w_{1 t_1})|
 \frac{1}{\eps^{k-n}} dz_{n+1 t_{n+1}}\ldots dz_{k t_k}\leq
 \left(\frac{C}{\delta}\right)^{k-n}.
\end{equation*}}
$\bullet$ $v_{m+1}^\eps,\ldots,v_n^\eps\in\Theta_2$.
 \jump Let $\eps^\beta\leq L\leq \delta$, and define the annulus,
$A(v_j^\eps,L)=\{v\in G^\eps|L\leq |v-v_j^\eps|\leq L+\eps\}$. By
definition of the paths $\gamma_j^\eps$, if $v_j^\eps\in
A(v_{j+1}^\eps,L)$, we have $\lim_{\eps\rightarrow
0}\frac{|b_j-w_{j+1}|}{\eps}=\infty$ (since $\beta<1$). Using
Theorem \ref{thm15} yields
\begin{equation*}
\frac{1}{\eps}|K^{-1}(b_j,w_{j+1})|=O\left(\frac{1}{|b_{j}-w_{j+1}|}\right)\leq
\frac{1}{\min_{\{w_{j+1}\in\gamma_{j+1}^\eps\}}|b_j-w_{j+1}|}=
\frac{C}{x_j+CL}\sin\theta_j,
\end{equation*}
where $x_j$ is the distance from $v_j^\eps$ to $b_j$. Let us
replace $\sin\theta_j$ by $C$. Hence, if for $m+1\leq j\leq n$,
$v_j^\eps\in A(v_{j+1},L_j)$, we obtain {\scriptsize
\begin{eqnarray*}
\sum_{|t_{m+1}|\leq\delta,\ldots,|t_n|\leq\delta}&&
 |K^{-1}(b_{m+1 t_{m+1}},w_{m+2 t_{m+2}})\ldots K^{-1}(b_{n t_n},w_{n+1 t_{n+1}})|
 \frac{1}{\eps^{n-m}} dz_{m+1 t_{m+1}}\ldots dz_{n t_n}\leq\\
&\leq&\sum_{|t_{m+1}|\leq\delta,\ldots,|t_n|\leq\delta}
 \frac{C}{x_{m+1 t_{m+1}}+C L_{m+1}}\ldots\frac{C}{x_{n t_n}+C
 L_n}dx_{m+1 t_{m+1}}\ldots dx_{n t_n},\\
&\sim& C\prod_{j=m+1}^n \log\left(\frac{\delta +C L_j}{C
L_j}\right).
\end{eqnarray*}}

\noindent Let $\Xi$ be the sum (\ref{534}) over vertices
$v_1^\eps,\ldots,v_k^\eps\in\Theta_1\cap\Theta_2\cap\Theta_3$.
Denote by $M=\sup_{v\in\RR^2}|\vphi(v)|$. Then, $\Xi\leq \Xi_1
\Xi_2 \Xi_3$, where {\scriptsize
\begin{eqnarray*}
\Xi_1&=&CM^m
\left[\sum_{v_1^\eps,\ldots,v_m^\eps\in\Theta_1}\eps^{2m}\left(\frac{\delta}{\eps}\right)^m\right],\\
\Xi_2&=&
 \left[\sum_{v_{n+1}^\eps,\ldots,v_k^\eps\in\Theta_2}\eps^{2(k-n)}|\vphi(v_{n+1}^\eps)|\ldots
 |\vphi(v_k^\eps)|
 \left(\frac{1}{\delta}\right)^{k-n}\right],\\
\Xi_3&=&C M^{n-m}
 \left[\sum_{L_{m+1}=0}^\delta\ldots\sum_{L_n=0}^\delta\eps^{n-m}L_{m+1}\ldots L_n
 \log\left(\frac{\delta+CL_{m+1}}{L_{m+1}}\right)\ldots\log\left(\frac{\delta+CL_n}{L_n}\right)
 dL_{m+1}\ldots dL_n\right].
\end{eqnarray*}}
Moreover,
\begin{eqnarray*}
\Xi_1&\leq& CM^m\delta^m \eps^{m(2\beta-1)},\\
\Xi_2&\leq&\left(\frac{1}{\delta}\right)^{k-n}\sum_{v_{n+1}^\eps,\ldots,v_k^\eps\in\RR^2}
 \eps^{2(k-n)}|\vphi(v_{n+1}^\eps)|\ldots|\vphi(v_k^\eps)|\leq M^{k-n}\left(\frac{1}{\delta}\right)^{k-n},\\
\Xi_3&\sim& CM^{n-m}(\delta \log\delta)^{n-m}.
\end{eqnarray*}
Hence, $\Xi\leq CM^k \delta^{2n-k}\eps^{m(2\beta-1)}$. Let us take
$\beta=2/3<1$, then $m(2\beta-1)>0$. If $2n-k\geq 1$, then
$\Xi=o(1)$. If $2n-k\leq 0$, take $\eps\leq
\delta^{\frac{k-2n+1}{m(2\beta-1)}}$, and $\Xi=o(1)$.

\begin{rem}\label{rem510}
Let $G^*$ be a bipartite isoradial graph, and let $K^{-1}$ be the
corresponding inverse Dirac operator, then for every black vertex
$b$ and every white vertex $w$ of $G^*$, we have
$$
\left|K^{-1}(b,w)\right| \leq \cal{C}
$$
for some constant $\C$ which only depends on the graph $G^*$.
\end{rem}
\begin{proof}
By theorem $4.2$ of \cite{Kenyon3}, $K^{-1}$ is given by
$$
K^{-1}(b,w)=\frac{1}{4 \pi^2 i}\int_{C}f_{wb}(z) \log z \,dz,
$$
where $C$ is a closed contour surrounding cclw the part of the circle
$\{e^{i\theta}|\theta\in[\theta_0-\pi+\Delta,\theta_0+\pi-\Delta]\}$,
which contains all the poles of $f_{wb}$, and with the origin in its exterior.
Without
loss of generality suppose $\theta_0=0$. As in \cite{Kenyon3}, let
us homotope the curve $C$ to the curve from $-\infty$ to the
origin and back to $-\infty$ along the two sides of the negative
real axis. On the two sides of this ray, $\log z$ differs by $2\pi
i$, hence
$$
K^{-1}(b,w)=\frac{1}{2 \pi}\int_{-\infty}^0 f_{wb}(t) \,dt, \mbox{
where } f_{wb}(t)=\frac{1}{(t - e^{i \theta_1})(t - e^{i
\theta_2})} \prod_{j=1}^k \frac{(t - e^{i \alpha_j})}{(t - e^{i
\beta_j})}.
$$
Refer to \cite{Kenyon3} for the choice of path from $w$ to $b$,
that is for the definition of the angles $\theta_1,
\theta_2,\alpha_j,\beta_j$. These angles have the property that
for all $j$, $\cos \alpha_j \leq \cos \beta_j$, so
$$
\left|\frac{t-e^{i \theta_j}}{t-e^{i \beta_j}}\right| \leq 1.
$$
For every $e^{i \theta} \in \{e^{i\theta}\,|\, \theta \in
[\theta_0-\pi+\Delta,\theta_0+\pi-\Delta]\}$, and for every $t
<0$, we have $|t-e^{i \theta}|^2 \geq
|t-e^{i(\pi+\Delta)}|^2=|t+e^{i\Delta}|^2$. Moreover for every $t
\in \RR$, we have $|t + e^{i \Delta}|^2 \geq \sin^2 \Delta$, and
$|t + e^{i \Delta}|^2 \geq (t+1)^2$, thus
$$
\int_{-\infty}^0|f_{wb}(t)| dt \leq \int_{-\infty}^{-2}
\frac{1}{(t+1)^2} dt + \int_{-2}^0 \frac{1}{\sin^2 \Delta} dt = 1
+ \frac{2}{\sin^2 \Delta}.
$$
Hence $ \left|K^{-1}(b,w)\right|\leq \C, \mbox{ where }
\C=\frac{1}{2 \pi}\left( 1+\frac{2}{\sin^2 \Delta}\right). $
\end{proof}
\end{proof}
We end the proof of Proposition \ref{prop53} with the following
\begin{lem}\
\begin{enumerate}
 \item[$1.$] $\displaystyle \int_{\RR^2}\ldots\int_{\RR^2}\vphi(v_1)\ldots \vphi(v_k)
\lim_{\eps \rightarrow 0}\EE[(h(v_1^\eps)-h(u_1^\eps))\ldots
(h(v_k^\eps)-h(u_k^\eps))]dz_1\ldots dz_k=$ {\flushright
$\qquad\qquad\qquad\qquad\qquad\qquad$ $ $ $=\left\{
\begin{array}{ll}
\displaystyle 0 & \mbox{ when $k$ is odd},\\
\displaystyle (k-1)!!\frac{1}{\pi^{k/2}}G(\vphi,\vphi)^{k/2} &
\mbox{ when $k$ is even}.
\end{array}
\right.$}
 \item[$2.$] $\lim_{\eps \rightarrow 0}\EE[H_{u_1}^\eps
\vphi \ldots H_{u_k}^\eps \vphi]=\lim_{\eps \rightarrow 0} \EE
[(H^\eps \vphi)^k]$.
\end{enumerate}
\end{lem}
\begin{proof}
 $1.$ is deduced from the formula of Proposition \ref{prop52}, and
from the fact that $\vphi$ is a mean $0$ function.
 $2.$ is a consequence of the fact that $\vphi$ is a mean $0$ function, and
of estimates of the kind of those of Lemma \ref{lem59}.
\end{proof}
\end{proof}

\subsection{Remark} \label{subsec45}

Note that the double periodicity assumption of the graph $G^*$ is only
required in Lemma \ref{lem54}, where we implicitly use the expression of
Theorem \ref{thm11} for the Gibbs measure $\mu$, as a function of the Dirac
operator $K$, and its inverse $K^{-1}$. By Remark \ref{rem1}, this
assumption can be released for the measure $\mu$, in the case where the
graph $G$ is a lozenge-with-diagonals tiling. Hence, Theorem
\ref{thm1} remains valid when the graph $G$ is any
lozenge-with-diagonals tiling of the plane, periodic or not.

\section{Proof of Corollary \ref{thm2}}\label{sec5}

We place ourselves in the context of Corollary \ref{thm2}: $\Q$ is the
set of triangular quadri-tilings, and assume quadri-tiles are assigned the
critical weight function. Recall the following notations: $\TT^\eps$
is the equilateral triangular lattice whose edge-lengths have been
multiplied by $\eps$, $\PP$ is the Gibbs measure on $\Q$ of Section
\ref{subsec14}; for $i=1,2$,
$H_i^\eps\vphi=\eps^2\sum_{v\in V(\TT^\eps)}\frac{\sqrt{3}}{2}\vphi(v)h_i^\eps(v)$, and
$F_1,F_2$ are Gaussian free fields of the plane.\jump
In order to
prove weak convergence in distribution of the height functions $h_1^\eps$
and $h_2^\eps$ to two independent Gaussian free fields $F_1$ and $F_2$, it suffices to show
that, $\forall \vphi\in\K$:
\begin{equation*}
\lim_{\eps\rightarrow 0}\EE[(H_1^\eps\vphi)^k (H_2^\eps\vphi)^m]=
\EE[(F_1\vphi)^k]\EE[(F_2\vphi)^m].
\end{equation*}
The key point is to obtain the analog of the moment formula of
Proposition \ref{prop51}. The rest of the proof goes through in 
the same way, and since notations are quite heavy, we do not repeat it here.\jump
The idea to obtain the moment formula is the following. Recall that triangular 
quadri-tilings correspond to two superposed
dimer models, the first on lozenge-with-diagonals tilings of $\LL$ and
the second on the equilateral triangular lattice $\TT$. Recall also
that both lozenge-with-diagonals tilings and $\TT$ are isoradial
graphs. Hence, we start by applying Proposition \ref{prop51} in the case
where the graph $G$ is a lozenge-with-diagonals tiling $L\in\LL$. Then
in Lemma \ref{lem62}, we prove some uniformity of convergence for every
$L\in \LL$, and we conclude the proof by using Proposition \ref{prop51}
in the case where the graph $G$ is the equilateral triangular lattice
$\TT$.\jump
Let $u_1,\ldots,u_k,v_1,\ldots,v_k,\us_1,\ldots,\us_m,\vs_1,\ldots,\vs_m$
be distinct points of $\RR^2$, and let
$\gamma_1,\ldots,\gamma_k$, $\gamma_1',\ldots,\gamma_m'$ be pairwise
disjoint paths such that $\gamma_j$ (resp. $\gamma_j'$) runs from
$u_j$ to $v_j$ (resp. from $\us_j$ to $\vs_j$). Define
\begin{equation*}
\GG(u_1,v_1,\ldots,u_k,v_k)=
\frac{(-i)^k}{(2\pi)^k}\sum_{\eps=0,1}(-1)^{k\eps}\left(
\int_{\gamma_1}\ldots\int_{\gamma_k} \det_{
\begin{array}{l}
\scriptstyle i,j \in [1,k]\\
\scriptstyle i \neq j
\end{array}} \left(
\frac{1}{z_i^\eps-z_j^\eps} \right) dz_1^\eps\ldots
dz_k^\eps\right).
\end{equation*}
Similarly, define $\GG(\us_1,\vs_1,\ldots,\us_m,\vs_m)$. Let
$u_j^\eps,v_j^\eps,\us_j^\eps,\vs_j^\eps$ be vertices of $\TT^\eps$
lying within $O(\eps)$ of $u_j,v_j,\us_j,\vs_j$ respectively.
\begin{lem}\label{lem61}
\begin{equation*}
\begin{split}
\lim_{\eps\rightarrow 0}
\EE[(h_1(v_1^\eps)-h_1&(u_1^\eps))\ldots(h_1(v_k^\eps)-h_1(u_k^\eps))
(h_2(\vs_1^\eps)-h_2(\us_1^\eps))\ldots(h_2(\vs_m^\eps)-h_2(\us_m^\eps))]=\\
&=\GG(u_1,v_1,\ldots,u_k,v_k)\GG(\us_1,\vs_1,\ldots,\us_m,\vs_m).
\end{split}
\end{equation*}
\end{lem}
\begin{proof}
By definition of the measure $\PP$, we have:\vspace{0.2cm}\\
$\EE[(h_1(v_1^\eps)-h_1(u_1^\eps))\ldots(h_1(v_k^\eps)-h_1(u_k^\eps))
(h_2(\vs_1^\eps)-h_2(\us_1^\eps))\ldots(h_2(\vs_m^\eps)-h_2(\us_m^\eps))]=$
{\scriptsize
\begin{eqnarray*}
&=& \sum_{\Ls^*\in\M(\TT^*)}
 \EE_{\mu^L}[(h_1(v_1^\eps)-h_1(u_1^\eps))\ldots(h_1(v_k^\eps)-h_1(u_k^\eps))
 (h_2(\vs_1^\eps)-h_2(\us_1^\eps))\ldots(h_2(\vs_m^\eps)-h_2(\us_m^\eps))
 ]d\mu^\TT(\Ls^*),\\
&=& \sum_{\Ls^*\in\M(\TT^*)}
 (h_2(\vs_1^\eps)-h_2(\us_1^\eps))\ldots(h_2(\vs_m^\eps)-h_2(\us_m^\eps))
 \EE_{\mu^L}[(h_1(v_1^\eps)-h_1(u_1^\eps))\ldots(h_1(v_k^\eps)-h_1(u_k^\eps))]d\mu^\TT(\Ls^*).
\end{eqnarray*}}
Using Section \ref{subsec45}, we can use Proposition
\ref{prop51} in the case where $G$ is any lozenge-with-diagonals tiling $L\in\LL$
(periodic or not), hence for every $L\in\LL$, we have:
\begin{equation}\label{61}
\lim_{\eps\rightarrow 0}
\EE_{\mu^L}[(h_1(v_1^\eps)-h_1(u_1^\eps))\ldots(h_1(v_k^\eps)-h_1(u_k^\eps))]=
\GG(u_1,v_1,\ldots,u_k,v_k).
\end{equation}
Note that the right hand side is independent of $L$, hence to obtain Lemma
\ref{lem61}, we need to prove that convergence in (\ref{61}) is
uniform in $L$, see Lemma \ref{lem62} below. Indeed, assuming this is the case, for
$\eps$ small, we can write:\vspace{0.2cm}\\
{\scriptsize
$
\displaystyle
\sum_{\Ls^*\in\M(\TT^*)}
 (h_2(\vs_1^\eps)-h_2(\us_1^\eps))\ldots(h_2(\vs_m^\eps)-h_2(\us_m^\eps))
 \EE_{\mu^L}[(h_1(v_1^\eps)-h_1(u_1^\eps))\ldots(h_1(v_k^\eps)-h_1(u_k^\eps))]d\mu^\TT(\Ls^*)
=$}
{\small
\begin{eqnarray*}
&=& \sum_{\Ls^*\in\M(\TT^*)}
 (h_2(\vs_1^\eps)-h_2(\us_1^\eps))\ldots(h_2(\vs_m^\eps)-h_2(\us_m^\eps)) 
 (\GG(u_1,v_1,\ldots,u_k,v_k)+O(\eps))d\mu^\TT(\Ls^*),\\ 
&=& (\GG(u_1,v_1,\ldots,u_k,v_k)+O(\eps)) 
 \sum_{\Ls^*\in\M(\TT^*)}
 (h_2(\vs_1^\eps)-h_2(\us_1^\eps))\ldots(h_2(\vs_m^\eps)-h_2(\us_m^\eps))d\mu^\TT(\Ls^*),\\
&=& (\GG(u_1,v_1,\ldots,u_k,v_k)+O(\eps))
 \,\EE_{\mu^\TT}[(h_2(\vs_1^\eps)-h_2(\us_1^\eps))\ldots
 h_2(\vs_m^\eps)-h_2(\us_m^\eps)],\vspace{0.1cm}\\
&=& (\GG(u_1,v_1,\ldots,u_k,v_k)+O(\eps))\,(\GG(\us_1,\vs_1,\ldots,\us_m,\vs_m)+O(\eps)),
\end{eqnarray*}}
where the last line is obtained by using Proposition \ref{prop51} for the
graph $\TT$.
\begin{lem}\label{lem62}
When $\eps$ is small, and for every lozenge-with-diagonals tiling
$L\in \LL$,
\begin{equation*}
\EE_{\mu^L}[(h_1(v_1^\eps)-h_1(u_1^\eps))\ldots(h_1(v_k^\eps)-h_1(u_k^\eps))]=
\GG(u_1,v_1,\ldots,u_k,v_k)+O(\eps),
\end{equation*}
where $O(\eps)$ is independent of $L$.
\end{lem}
\begin{proof}
Let us look at the proof of Proposition \ref{prop51} in the case of a
lozenge-with-diagonals tiling $L\in\LL$, and denote by $K_L$ the Dirac
operator indexed by vertices of $L^*$. Then Lemma \ref{lem62} is
proved if we show that $O(\eps)$ in Lemma \ref{lem56} is independent
of $L\in\LL$. Looking at the proof of Lemma \ref{lem56}, we see that
$O(\eps)$ comes from the error term in the asymptotic formula for the inverse
Dirac operator of Theorem \ref{thm15}, \cite{Kenyon3}:
\begin{equation*}
K_L^{-1}(b_i,w_j)=\eps\left(\frac{1}{2 \pi} (F_0(b_i,w_j)+ f_{w_j
b_i}(0)F_1(b_i,w_j))+O(\frac{\eps}{|b-w|^2})\right).
\end{equation*}
In \cite{Kenyon3}, the error term is computed
explicitly. Looking at the explicit formula, and using the
regularity of the graphs $L$, we show that
$O(\frac{\eps}{|b-w|^2})=\frac{C_1\eps}{|b-w|^2}$, where $C_1$ is
independent of $L$. Moreover, by assumption the paths
$\gamma_1,\ldots,\gamma_k$ are disjoint, so that we define:
\begin{equation*}
C_2=\inf_{i\neq j}\inf_{\{b\in \gamma_i,w\in\gamma_j\}}|b-w|>0, 
\end{equation*}
which is independent of $L$. Hence $O(\frac{\eps}{|b-w|^2})=\frac{C_1\eps}{C_2^2}$.
\end{proof}
\end{proof}

\bibliographystyle{alpha}

\end{document}